\documentclass[12pt,oneside,a4paper,leqno]{article}
\usepackage{amsfonts,amssymb}
\usepackage[centertags]{amsmath}
\usepackage[all]{xy}\xyoption{all}
\usepackage{amsthm}
\usepackage{stmaryrd}
\newcounter{lemma}
\swapnumbers
\theoremstyle{plain}
\newtheorem{lemma}[equation]{Lemma}\numberwithin{lemma}{section}
\newtheorem{theo}[equation]{Theorem}
\newtheorem{propo}[equation]{Proposition}
\newtheorem{coro}[equation]{Corollary}
\theoremstyle{definition}
\newtheorem{defi}[equation]{Definition}
\newtheorem{example}[equation]{Example}
\newtheorem{remark}[equation]{Remark}
\theoremstyle{remark}

\numberwithin{equation}{section}

\makeatletter
\def\tagform@#1{\maketag@@@{\ignorespaces#1\unskip\@@italiccorr}}
\renewcommand{\theequation}{(\arabic{section}.\arabic{equation}\/)}

\makeatother
\newcommand{\ldiag}[2]{%
\begin{equation}\label{#1}\begin{aligned}\xymatrix{#2}\end{aligned}\end{equation}%
}
\newcommand{\diag}[1]{%
\begin{equation}\begin{aligned}\xymatrix{#1}\end{aligned}\end{equation}%
}
\newcommand{\ndiag}[1]{%
\begin{equation*}\begin{aligned}\xymatrix{#1}\end{aligned}\end{equation*}%
}
\newdir{ >}{!/8pt/@{ }*@{>}}
\newdir{< }{!/-8pt/@{ }*@{<}}
\newcommand{\hh}{\Phi}
\newcommand{\rel}{\  \mathrm{rel}\;}

\newcommand{\minus}{\smallsetminus\/}

\newcommand{\R}{\mathbb{R}}

\newcommand{\V}{{\bf\mathfrak{{F}}}}

\newcommand{\Ab}{\mathbf{Ab}}

\newcommand{\Vtop}{\mbox{$\V$-${{{\mathbf{Top}}}}$}}

\newcommand{\homotopic}{\sim}
\renewcommand{\top}{\mathbf{Top}} 

\newcommand{\T}{\mathbb{T}}

\newcommand{\nota}[1]{}%\marginpar{\fbox{\large\bf\footnote{#1}}}}

\newcommand{\op}{\mathrm{op}}

\newcommand{\from}{\colon}

\renewcommand{\equiv}{\simeq}
\newcommand{\pr}{\mathrm{pr}}

\newcommand{\ff}{F} % fibre functor
\newcommand{\grd}{\mathbf{Grd}}
\newcommand{\gpd}{\mathbf{Grd}}
\newcommand{\FAgrd}{(\V,A)\mbox{-}\mathbf{Grd}}
\newcommand{\FAset}{(\V,A)\mbox{-}\mathbf{Set}}
\newcommand{\ze}{\mathbb{Z}}
\newcommand{\Mod}{\mathbf{Mod}}
\newcommand{\Mor}{\mathrm{Mor}}
\newcommand{\CC}{\mathbf{C}}

\newcommand{\Set}{\mathbf{Set}}
\newcommand{\model}{\mathbf{model}}
\newcommand{\coef}{\mathbf{Coef}}

\newcommand{\tto}{\Longrightarrow}
\newcommand{\PP}{\mathbf{P}}
\newcommand{\expansion}{\nearrow}
\newcommand{\collapse}{\searrow}
\newcommand{\deformation}{\curvearrowright}
\newcommand{\Wh}{\mathbf{Wh}}
\newcommand{\torsion}{\tau}

\newcommand{\obst}{\mathcal{O}}
\newcommand{\TT}{\mathbf{T}}
\newcommand{\Kiso}{K^\mathrm{iso}_1}
\newcommand{\AAA}{\mathbf{A}}

\newcommand{\ringoids}{\mathbf{Ringoids}}
\newcommand{\image}{\mathrm{Image}}

\newcommand{\PiV}{\llbracket\V\rrbracket}
\begin{document}
\pagenumbering{arabic}

\title{%
Homotopy and homology of fibred spaces
}

\author{Hans-Joachim Baues and Davide L.~Ferrario}
% \email{baues@mpim-bonn.mpg.de}
% \address{
% Max-Planck-Institut f\"ur Mathematik\\
% Vivatsgasse, 7 - 53111 Bonn (DE) 
% }
%  \email{ferrario@mpim-bonn.mpg.de}

% \author{Davide L.~Ferrario}
 % \email{ferrario@mpim-bonn.mpg.de}

%\address{
% Max-Planck-Institut f\"ur Mathematik\\
% Vivatsgasse, 7 - 53111 Bonn (DE) 
% }

% \date{\today}
% \centerline{\bgroup \large \today \ --- \now \egroup}
\date{}

\maketitle

\begin{abstract}
We study fibred spaces with fibres in a structure category $\V$
and we show that cellular approximation, Blakers--Massey theorem,
Whitehead theorems, obstruction theory,
Hurewicz homomorphism, Wall finiteness obstruction, and Whitehead
torsion theorem hold for fibred spaces. For this we introduce the 
cohomology of fibred spaces.
\end{abstract}

\emph{AMS SC}:
55R55 (Fiberings with singularities);
54H15 (Transformation groups and semigroups).
55R65 (Generalizations of fiber spaces and bundles);
55R70 (Fiberwise topology);

%======================================================================
\section{Introduction}
A fibred space $X$ is a map 
\[
p\from X \to \bar X
\]
for which fibres $p^{-1}x$, $x\in \bar X$, are objects in a structure
category $\V$. For example a $G$-space $X$ for a topological
group $G$ which acts properly on $X$ yields the fibred space
\[
p\from  X \to X/G
\]
where $X/G$ is the orbit space. Here fibres are objects in the orbit
category $\V_G$ consisting of $G$-orbits $G/H$ and $G$-equivariant maps.
Moreover fibre bundles \cite{huse,steenrod}
and stratified fibre bundles
\cite{qq} are examples of fibred spaces.
In this paper we study the homotopy theory of fibred spaces. 
In particular we show that fundamental results like Whitehead theorems,
cohomology, obstructions to the extension of maps,
Wall-finiteness obstruction and Whitehead torsion are available
for fibred spaces.
These results form the bulk of the book \cite{lueck} in the case 
of $G$-spaces. As it turns out homotopy theory of $G$-spaces 
to a large extent is just a specialization of the homotopy
theory of fibred spaces.

The theory of $\V$-fibred spaces yields a new flexibility since the structure
category $\V$ is a parameter of the theory. For example all constructions and
obstructions are natural in $\V$ and such naturality is useful 
similarly as in the 
case of
naturality of coefficients in cohomology.  By comparing the orbit category
$\V_G$ with other structure categories $\V$ one also obtains 
applications to
the theory of transformation groups. On the other hand our theory can be
applied to stratified bundles like stratified vector bundles which appear as
tangent bundles of stratified manifolds, see \cite{kreck,sm,qq}.
Obstructions to non trivial sections of stratified vector bundles
can be studied by use of the obstruction theory of fibred spaces (such problems
were treated in \cite{schwartz}).

A self-contained proof of the results in this paper is highly elaborate and
requires a rewriting of the book \cite{lueck} (more than 400 pages) for fibred
spaces. In this paper, however, we prove the results by use of the axiomatic
approach of \cite{baues}. For this we only need to show that the
Blakers--Massey theorem holds in the category of fibred spaces.  This implies
that we can use the methods in \cite{baues} leading, in particular, to a new
construction of the cellular chain complex $C_*(X,A)$ which is simpler than the
construction of Bredon \cite{bredon,brocker,tomdieck,lueck}. 
%\cite{brocker}, \cite{tomdieck}, \cite{lueck}. 

Our homotopy theory of fibred spaces with fibres in the (topological enriched)
category
$\V$ is closely related to the theory of continuous $\V^\op$-diagrams of spaces
(compare with the principal diagram in \cite{qq}). 
In case $\V$ is discrete such diagrams are treated in chapter A of
\cite{baues}. It is easy to transform the results in this paper to 
the theory of continuous diagrams by use of \cite{baues}.

The methods of \cite{baues}, in fact, yield further results on fibred spaces
not treated in this paper. We leave it to the reader to study the ``model
lifting property'' of the chain complex $C_*(X,A)$, the ``obstruction for the
realizability'' of a chain complex, the ``tower of categories'' for the
classification of homotopy classes of maps between fibred spaces and the
associated spectral sequences.

On the other hand we refer to the recent book of Crabb and James
\cite{crabbjames} for a different perspective on 
spaces of fibres 
and for a detailed account of the areas
of 
fibrewise topology and fibrewise homotopy theory.

\section{Fibre families}
\label{section:fibrespaces}
Let $\V$ be a small category  
together with a faithful 
functor $\ff\from \V \to  \top$ to the category 
$\top$ of compactly generated Hausdorff spaces. 
Then $\V$ is 
termed 
a \emph{structure category}  
and $\ff$ is the \emph{fibre functor} on $\V$.
In many examples the functor $\ff$ is actually
the inclusion of a subcategory $\V$ of $\top$ 
so that in this case we need not to mention the fibre
functor $\ff$.

A 
\emph{fibre family}
with fibres in $\V$ (or a
$(\V,\ff)$-family)
is 
a topological space $X$, termed \emph{total space},
together with a map $p_X\from X \to \bar X$,
termed \emph{projection}
to the \emph{base space} $\bar X$, 
and for every $b\in \bar X$ a selected
homeomorphism $\Phi_b\from p_X^{-1}b \approx \ff X_b$
where $X_b$ is an object in $\V$, called \emph{fibre},
depending on $b\in \bar X$.
The homeomorphism $\Phi_b$ is termed \emph{chart} at $b$.
The family $(p_X\from X \to \bar X, X_b, \Phi_b, b\in \bar X)$
is denoted simply by $X$.
A fibre family is also termed a \emph{fibred space} with fibres
in $\V$.

Given two $\V$-families $X$ and $Y$ a \emph{$\V$-map} from $X$ to $Y$
is a pair of maps $(f,\bar f)$ such that the following diagram
\ndiag{% 
X \ar[r]^f
\ar[d]_{p_X} & Y \ar[d]^{p_Y} \\
\bar X \ar[r]^{\bar f} & \bar Y,
}%
commutes, and such that for every $b\in \bar X$ 
the composition given by the dotted arrow of the diagram
\[
\xymatrix@C+24pt{% 
p_X^{-1}(b)\ar[r]^{f|p^{-1}(b)}
&  p_Y^{-1}(\bar f b) \ar[d]^{\Phi_{\bar f b}} \\
\ff X_b 
\ar[u]^{\Phi_b^{-1}} 
\ar@{.>}[r]^{} & \ff Y_{\bar f b}
}%
\]
is a morphism 
in the image of the functor $\ff$. That is,
there exists a morphism $\phi\from X_b \to Y_{\bar f b}$
in $\V$
such that the dotted arrow is equal to $\ff(\phi)$.
We will often denote $\ff X_b$ by $X_b$
and it will be clear from the context whether $X_b$ 
denotes an object in $\V$ or a space in $\top$ given by
the functor $\ff$.
If a $\V$-map $f=(f,\bar f)$ is a $\V$-isomorphism
then $f$ and $\bar f$ are homeomorphisms but the converse need not be true.

If $X$ is a $\V$-family  and $Z$ is a 
topological space, then the product $X\times Z$ in $\top$ is a $\V$-family
with projection 
$p_{X\times Z} = p\times 1_Z\colon X\times Z \to \bar X \times Z$.
The fibre over  a point $(b,z)\in \bar X\times Z$ is equal to 
$p_{X}^{-1}(b)\times \{z\}$; using 
the chart $\Phi_b:p_{X}^{-1}(b) \to X_b$
the chart $\Phi_{(b,z)}$ is defined by $(x,z) \mapsto \Phi_b(x) \in X_b$
where of course we set $(X\times Z)_{(b,z)} = X_b$. 
In particular, by taking $Z=I$ the unit interval
we obtain the cylinder object $IX=X\times I$
% $\Vtop$ has the cylinder functor $-\times I\colon \Vtop \to \Vtop$,
and therefore the notion of \emph{homotopy}: two $\V$-maps
$f_0,f_1\colon X \to Y$ 
are $\V$-homotopic (in symbols $f_0\homotopic f_1$)
if there is a $\V$-map $F\colon IX \to Y$ 
such that $f_0=Fi_0$ and $f_1=Fi_1$. Here $i_0$ and $i_1$ 
are the inclusions $X \to X\times I$ at the levels $0$ and 
$1$ respectively.

Let $\Vtop$ be the category consisting 
of $\V$-families $p_X: X \to \bar X$ in $\top$ and $\V$-maps.
Homotopy of $\V$-maps yields a natural equivalence relation $\sim$
on $\Vtop$ so that the homotopy category
$(\Vtop)/{_\homotopic}$ is defined.

\begin{defi}
If $V$ is a fibre in $\V$ and $\bar X$ is a space in $\top$,
then the projection onto the first 
factor 
$p_1\from X=\bar X\times V \to \bar X$ 
yields the \emph{product family} 
with fibre $V$;
the charts $\Phi_b\from \{b\}\times V \to V=X_b$
are given by projection and $X_b=V$ for all $b\in \bar X$.
If $X$ is a $\V$-family $\V$-isomorphic to a product
family then $X$ is said to be a \emph{trivial $\V$-bundle}.
In general a \emph{$\V$-bundle} is a locally trivial family of 
fibres, i.e.  a family $X$ over $\bar X$ such that 
every $b\in \bar X$ admits a neighborhood $U$ for which $X|U$ is trivial.
Here 
$X|U$
is the \emph{restriction}
of the family $X$
defined by $U\subset \bar X$.
\end{defi}

% mapping space

Given a family $Y$ with projection $p_Y\colon Y \to \bar Y$ and 
a map $\bar f\colon \bar X \to \bar Y$, the pull-back $X=\bar f^*Y$
is the total space of a family of fibres given by the vertical
dotted arrow of the following pull-back diagram in $\top$.
\diag{%
X =
\bar f^* Y 
\ar@{.>}[r]
\ar@{.>}[d]
&  Y
\ar[d]^{p}
\\
\bar X 
\ar[r]^{\bar f}
& \bar Y.\\
}%
The charts are defined as follows:  For every $b\in \bar X$ 
let $X_b= Y_{\bar fb}$, and $\Phi_b: p_X^{-1}(b) \to X_b$ 
the composition 
$p_X^{-1}(b)\to p_Y^{-1}(\bar fb) \approx Y_{\bar fb}=X_b$
where the map $p_X^{-1}(b)\to p_Y^{-1}(\bar fb)$
is a homeomorphism since $X$ is a pull-back.

A $\V$-map $i\from A \to Y$ is termed a \emph{closed
inclusion} if $\bar i\from  \bar A \to \bar Y$ is an inclusion, 
$\bar i \bar A$ is closed in $\bar Y$ and the following diagram is 
a pull-back:
\diag{%
i^*Y = A \ar[r]^{i} \ar[d]_{p_A} & Y \ar[d]^{p_Y}\\
\bar A \ar[r]^{\bar i} & Y \\
}%
Hence a closed inclusion
$i\from A \to Y$ induces homeomorphisms
on fibres.

The push-out construction can be extended to the category $\Vtop$,
provided the push-out is defined via  a closed inclusion,
see \cite{qq} (2.5). %\ref{qq_lemma:pushout}
\begin{lemma}\label{lemma:pushout}
Given $\V$-families $A$, $X$, $Y$ and $\V$-maps $f\from A \to X$,
$i\from A\to Y$ with $i$ a closed inclusion 
%todo shorten the map notation
the push-out diagram %$Z$ of $f$ and $i$ 
in $\Vtop$ 
\ndiag{%
A \ar[r]^f \ar@{->}[d]_{i} & X \ar@{->}[d]\\
Y \ar[r] & Z \\
}%
exists 
and 
$X\to Z$ is a closed inclusion. 
\end{lemma}

Let $A\subset Y$ be a closed inclusion then also
$A\times I \subset Y \times I$ is a closed inclusion and the
\emph{relative cylinder} $I_AY$ is defined by the push-out
\ndiag{%
A\times I \ar[r]^{\pr} \ar[d] & A \ar[d] \\
Y\times I \ar[r] & I_AY \\
}%
in $\Vtop$. A map $F\from I_AY\to X$ is termed 
a \emph{homotopy} $f_0\sim f_1 \rel A$ with $f_0 = Fi_0$ 
and $f_1=Fi_1$. Let $g\from A \to X$ given. Then
\begin{equation}
\label{eq:relative}
[Y,X]_\V^A = [Y,X]^g_\V = 
\left\{f\from Y\to X, f|A = g\right\}/_{\sim \rel A} 
\end{equation}
denotes the set of homotopy classes relative $A$ or under $A$.
We call a closed inclusion $(Y,Y')$ a \emph{pair} in $\Vtop$
and maps between pairs are defined as usual.
Moreover for a closed inclusion $A\subset Y'$ and a map $g\from A \to X'$ we obtain the set
$[(Y,Y'),(X,X')]_\V^A$ of homotopy classes of pair maps
under $A$.

\begin{remark}
All results in this paper remain true
if $\top$ is the category of all topological spaces
and if the fibre functor $\ff$ satisfies condition
$(*)$:
For
every object $V$ in $\V$
the space $\ff(V)$
is locally compact
and Hausdorff. 
In \cite{qq} we assumed the fibres to be also second countable,
in order to deal with metrizable spaces and avoid many unnecessary 
technicalities, which do not occur in homotopy theory
of $\V$-families.
\end{remark}

\section{$\V$-complexes}
\label{section:cwcomplexes}

For an object $V$ in $\V$ the 
family $\ff V\to *$ with base space a singleton is termed a
\emph{$\V$-point} also denoted by $V$.
A disjoint union  of $\V$-points
is called a \emph{$\V$-set}. 
This is a $\V$-family for which the base space has
the discrete topology.
Let $D^n$ be the unit disc in $\R^n$ and $S^{n-1}$ 
its boundary with base point $*\in S^{n-1}$. 
The complement $e^n=D^n\minus S^{n-1}$ is the \emph{open cell}
in $D^n$.
A $\V$-cell is a product family $V\times e^n\to e^n$ with $V\in \V$.

We say that a $\V$-family $X$ is obtained from a 
$\V$-family $A$ by \emph{attaching}  $n$-cells if 
a $\V$-set $Z$ together with a  $\V$-map $f$ is given,
such that the following  diagram
\ldiag{diag:cw}{%
Z\times S^{n-1} \ar[r]^f \ar[d] & A \ar[d] \\
Z \times D^{n} \ar[r]^{\hspace{-38pt} \hh} & X=A\cup_f(Z\times D^n) \\
}%
is a push-out in $\Vtop$.
The inclusion $Z\times S^{n-1} \to Z\times D^{n}$
is a closed inclusion, therefore the push-out exists
and the induced map $A \to X$  is a closed inclusion
and $X\minus A = Z\times e^n$ is a union of open $\V$-cells.
If $Z$ is a $\V$-point then we say that $X$ is obtained from $A$
by attaching a $\V$-cell and $\hh$ is the
\emph{characteristic map} of the $\V$-cell.

\begin{defi}
\label{defi:CWcomplex}
A \emph{relative $\V$-complex} $(X,A)$ is a family $X$ and a 
filtration
\begin{equation*}
A=X_{-1} \subset X_0 
\subset X_1 \subset \dots \subset X_n \subset X_{n+1} \subset \dots \subset X
\end{equation*}
of $\V$-families $X_n$, $n\geq -1$, such that 
% $X_0$ is the coproduct of $D$ and a $\V$-set $Z_0$, 
for every $n\geq 0$ the $\V$-family $X_n$ is obtained from 
$X_{n-1}$ by attaching $n$-cells and
\begin{equation*}
X = \lim_{n\geq 0} X_n.
\end{equation*}
\end{defi}

The spaces $X_n$ are termed \emph{$n$-skeleta} of $(X,A)$.
%$X_n$ is said the \emph{$n$-dimensional skeleton} or 
%equivalently $n$-skeleton. 
If $A$ is empty we call $X$ a 
% \emph{$\V$-stratified complex} or a 
\emph{$\V$-complex}.
Then $X$ is a union of $\V$-cells, that is, 
there are the $\V$-sets $Z_n$ of $n$-cells in $X\minus A$
such that $X\minus A$ 
is the union of $Z_n\times e^n$ for $n\geq 0$.

\begin{example}
\label{ex:orbit}
Let $G$ be a compact Lie group
and let $\V$ be the \emph{category of orbits} of $G$,
that is, 
$\V$ is the subcategory of $\top$ consisting of spaces
$G/H$, where $H$ is a closed subgroup of $G$,
and $G$-equivariant maps $G/H \to G/H'$.
%see \cite{tomdieck}.
Then each $G$-CW-complex  (see \cite{tomdieck})
is a $\V$-complex.
\end{example}

\begin{example}
Each $\V$-stratified fibre bundle with finitely many strata
is $\V$-homotopy equivalent to a $\V$-complex \cite{qq}.
\end{example}

A map $f\from (K,A) \to (X,A)$
under $A$ between relative $\V$-complexes is termed
\emph{cellular} if $f(K_n) \subset X_n$ for $n\geq 0$.

\begin{theo}[Cellular approximation]
\label{theo:cat}
Let $(K,A)$ and $(X,A)$ be relative $\V$-complexes; let $L$ be
a subcomplex of $K$ and $g\from K \to X$ a $\V$-map under
$A$ such that the restriction to $L\subset K$ is cellular.
Then there is a cellular map $f\from K\to X$ extending $g|L$ and a
homotopy $f\sim g \rel L$.
\end{theo}
\begin{proof}
% theorem in $\Vtop$ is a consequence of IV.5.8 of Baues \cite{baues}.
The proof is identical to the classical case of CW-complexes,
once lemma \ref{lemma:nconnected} is proved. Alternatively
one can use IV.5.8 Baues \cite{baues}
to obtain the result by use 
of the Blakers-Massey theorem \ref{theo:BMp} below.
\end{proof}

\begin{remark}
In \cite{qq} we have seen that $\V$-complexes are stratified bundles which 
generalize the classical notion of bundle. 
The bundle theorem of \cite{qq} shows:
%\label{theo:mainstra2} \label{theo:bundlestra} 
Let $\V$ be a structure category which is a groupoid.
%and assume that the fibre functor satisfies condition 
%$(*)$:
%\begin{quote}
%For
%every object $V$ in $\V$
%the space $\ff(V)$
%is locally compact, second-countable and Hausdorff.
%\end{quote}
Then a $\V$-complex $X \to \bar X$
is a $\V$-bundle over $\bar X$. Conversely,
each $\V$-bundle
$X\to \bar X$ over a CW-complex $\bar X$ is $\V$-isomorphic to 
a $\V$-complex.
\end{remark}

\begin{remark}
The definition of $\V$-complex follows the lines of 
\cite{baues0} (page 194). 
After Whitehead, Bredon defined $G$-complexes for $G$ finite \cite{bredon},
independently Matumoto \cite{matumoto} and Illmann \cite{illman}
extended the definition to 
the case $G$ a compact Lie group. See \cite{tomdieck}
for a general approach of $G$-CW complexes and \cite{lueck} for 
its consequences in $K$-theory of $G$-complexes.
A notion of fibrewise CW-complex was presented by James
in \cite{james2} (see also the decomposition
of Bredon $O(n)$-manifolds 
$W^{2n-1}_k$ there) in the framework of fibrewise topology \cite{james}.
While in the case of $G$-CW-complexes the structure category
is evidently the orbit category of $G$, the structure category
of a fibrewise complex \cite{james2} is the full subcategory of 
$\top$ containing the fibres and the fibres are assumed to be compact.
\end{remark}

\section{Homotopy groups}

A $\V$-family  $X$ in $\Vtop$  is \emph{pointed}
if a $\V$-map $*_V\from  V\to X$ is given,
where $V$ denotes the $\V$-point given by the object $V$ in $\V$.
Such a map is termed a 
\emph{base point} of $X$. A \emph{pointed pair} $(X,Y,*_V)$ of $\V$-families
is given by a closed inclusion $Y\subset X$ and a base point 
$*_V\from V \to Y$. As usual a base point is chosen in the sphere $S^{n-1}$
so that 
the pair $(D^n,S^{n-1},*)$ is pointed in $\top$.
For every $n\geq 0$ and $\V$-point $V$ let $S^n_V$ 
and $D^{n+1}_V$ denote the pointed $\V$-family
$S^n_V = V\times S^n$ over $S^n$ and $D^{n+1}_V = V \times D^{n+1}$ 
over $D^{n+1}$ with $S^n_V\subset D^{n+1}_V$.
\begin{defi}
For every pointed $\V$-family $(X,*_V)$ the \emph{homotopy group}
\begin{equation}
\pi_n^\V(X;*_V) = \left[ S^n_V, X \right]^V_\V
\end{equation}
is the set of homotopy classes relative $V$ of $\V$-maps 
$f\from S^n_V\to X$
as in the commutative diagram \ref{diag:pi} in $\Vtop$.
\ldiag{diag:pi}{%
&V\times * \ar@{=}[r] \ar@{ >->}[d] & V \ar[d]^{*_V}\\
S^n_V\ar@{=}[r]&V\times S^n \ar[r]^f & X\\
}%
Moreover the 
\emph{relative homotopy group}
\begin{equation}
\pi_{n+1}^\V(X,Y;*_V) = \left[ (D^{n+1}_V, S^n_V), (X,Y) \right]^V_\V
\end{equation}
is the set of homotopy classes relative $V$ of pair  
maps $f\from (D^{n+1}_V, S^n_V) \to (X,Y)$.
\end{defi}

The homotopy groups above can also be described as ordinary
homotopy groups of function spaces.
For a fibre family $X$ and $V$ in $\V$ let $X^V$ be the space 
of all $\V$-maps $V\to X$ with the compact-open topology in $\top$.
For example the faithful fibre functor $\ff$ yields
the identification
\[
W^V = \hom_\V(V,W)
\]
so that $\V$ is a topological enriched category.
Function spaces yield the functor
\begin{equation}
X^\circ \from \V^\op \to \top
\end{equation}
which carries $V$ to $X^V$ and $\alpha\from V\to W$ 
to the induced map $\alpha^*\from X^W \to X^V$ with
$\alpha^*(x) = x \circ \ff(\alpha)$.

We have canonical isomorphisms
\begin{equation}
\pi_n^\V(X;*_V) = \pi_n(X^V;*_V)
\end{equation}
\[
\pi_n^\V(X,Y;*_V) = \pi_n(X^V,Y^V;*_V)
\]
where the right hand side denotes the homotopy groups in $\top$.

We will prove in section \ref{sec:whitehead}
the following Whitehead theorem in $(\Vtop)^A$.
\begin{theo}
\label{theo:whitehead}
A map $f\from (X,A) \to (Y,A)$ under $A$ 
between relative $\V$-complexes
is a $\V$-homotopy equivalence under $A$ if and only if  $f$ 
induces isomorphisms
\[
f_*\from \xymatrix@1{\pi_n^\V(X;*_V) \ar[r]^{\equiv} & \pi_n^\V(Y;f*_V)}
\]
for all $n\geq 0$, basepoints $*_V\from V\to X$, and objects
$V$ in $\V$.
\end{theo}

\begin{coro}
\label{coro:james}
Let $f$ be a   $\V$-map $f\from (X,A) \to (Y,A)$ under $A$ 
between relative $\V$-complexes.
Then the following statements are equivalent.
\begin{enumerate}
\item \label{c1}
The $\V$-map $f$ is a $\V$-homotopy equivalence under $A$;
\item \label{c2}
For every $V$ in $\V$ the induced map $f^V\from X^V \to Y^V$ 
is a weak homotopy equivalence in $\top$.
\end{enumerate}
\end{coro}
\begin{proof}
It is clear that \emph{(i)} $\implies$ \emph{(ii)}. %$\implies$ \emph{(iii)}.
But if $f^V$ is a weak equivalence in $\top$ for every
$V$, then by \ref{theo:whitehead} \emph{(i)} holds.
\end{proof}

\begin{remark}
It is clear that  in \ref{coro:james} we can replace 
\emph{(ii)} 
with the following:
For every $V$ in $\V$ the induced map $f^V\from X^V \to Y^V$ 
is a homotopy equivalence in $\top$. 
Thus we obtain a generalization to $\V$-complexes of 
James--Segal theorem
\cite{jamessegal} for $G$-ANR's. See also \cite{hauschild}.
\end{remark}

\section{The Blakers-Massey theorem}
We say that a pair $(X,Y)$ is \emph{$n$-connected}
if 
for all basepoints $*_V\from V\to Y$, $V\in \V$ and 
$1\leq r\leq n$ the homotopy group
$\pi_r^\V(X,Y;*_V) = 0$ 
are trivial and $\pi_0^\V(Y;*_V) \to \pi_0^\V(X;*_V)$
is surjective.

\begin{lemma}\label{lemma:nconnected} 
For a relative $\V$-complex $(X,A)$ the pair
$(X,X_n)$ is $n$-con\-nec\-ted, $n\geq 0$.
\end{lemma}
\begin{proof}
Let $V\in \V$.
As in the classical case it suffices to prove that for an 
attachment $X=A\cup_\alpha(V\times D^n)$ as in \ref{diag:cw} the 
pair $(X,A)$ is $(n-1)$-connected. For this we use similar arguments
and notations as in the proof of 13.5 of \cite{gray}. Let 
$g\from v\times (B^r,S^{r-1}) \to (X,A)$ 
be a map which represents an element in 
$\pi_r^\V(X,A;*_V)$, with $r<n$.
Using the classical map $\chi\from (I^r,\partial I^r,J^{r-1}) 
\to (B^r,S^{r-1},*)$ we obtain the composite
$f'=g\chi$ which plays the role of $f$ 
in the proof of $13.5$ of \cite{gray}.
The map $f'$ induces the commutative diagram
\[
\xymatrix@C+24pt{%
V \times (I^r,\partial I^r, J^{r-1}) \ar[r]^{f'}  \ar[d] & 
(A\cup(V\times e^n),A,A_*) \ar[d] \\
(I^r,\partial I^r,J^{r-1}) \ar[r]^{f} & 
(\bar A \cup e^n, \bar A, *) \\
}%
\]
Here $*\in \bar A$ is determined by $*_V$ 
and $A_*\in \V$ is the fibre of $A$ over $*$.
As in 13.5 of \cite{gray} we obtain for $f$
the subspace $\bar U\subset I^r$
with $f(\bar U)\subset e^n$
and the homotopy $h_t\from \bar U \to e^n$ ($t\in [0,1]$)
relative $\partial \bar U$ with $h_0 = f|\bar U$.
Hence we obtain by $f'$ the commutative diagram
\[
\xymatrix@C+24pt{
V\times \bar U \ar[r]^{f'|(V\times \bar U)} \ar[d]^q 
& V\times e^n \ar[d]^{q} \\
\bar U \ar[r] ^{f|\bar U} & e^n \\
}%
\]
The map $f'|(V\times \bar U)$ 
is a $\V$-map which is determined by $(f|\bar U)_q$ 
and the 
coordinate $f''\from V\times \bar U \to V$.
We define a $\V$-homotopy 
\[
h_t'\from V\times \bar U \to V\times e^n
\]
by $h'_t = (h_tp,f'')$.
We have $h_0 = f'|(V\times \bar U)$ and $h'_t$ us a $\V$-homotopy
relative $V \times \partial \bar U$.
Hence we can define a $\V$-homotopy 
\[
H'_t\from V\times I^r \to A\cup (V\times e^n) \mbox{ \ \ with $H'_0=f'$}
\]
which induces $H_t\from I^r \to \bar A\cup e^n$
in the proof of 13.5 of \cite{gray}.
We obtain $H'_t$ by
\[
H'_t(x,u) =
\begin{cases}
h'_t(x,u) & \mbox{for $u\in \bar U$, $x\in V$} \\
f'(x,u) & \mbox{ for $u\in I^r\minus U$, $x\in V$}\\
\end{cases}
\]
Now we can choose $p\in e^n(1/2) \subset e^n$
with $p\not\in \mathrm{image(H_1)}$; see
13.5 \cite{gray}.
Hence the map $H'_1$ has a factorization
\[
H'_1\from V\times I^r \to A\cup(V\times (e^n\minus p) ) 
\subset A\cup (V\times e^n).
\]
Here the inclusion $A\to A\cup(V\times(e^n\minus p))$
is an $\V$-homotopy equivalence.
%by \ref{lemma:ghost} below,
%which is a consequence of the actioms of $I$-category.
This completes the proof.
\end{proof}

\begin{lemma}
\label{lemma:approximation}
Let $(X,A)$ be a relative $\V$-complex 
and let $(X,A)$ be $n$-connected,
with $n\geq 0$.
Then there exists a relative $\V$-complex
$(Y,A)$ with $Y_0=Y_1=\dots=Y_n = A$ 
and a $\V$-homotopy equivalence $Y\to X$ under $A$.
\end{lemma}
\begin{proof}
The map $Y\to X$ is obtained inductively
by attaching ``ball pairs'' and collapsing
the non-attached part of the boundary of the ball pair;
compare the proof of (6.14) in \cite{switzer}.
\end{proof}

\begin{theo}\label{theo:BMp}
Let $(X,A)$ a relative $\V$-complex with subcomplexes $X_1$,$X_2$,
such that $X=X_1\cup X_2$. Let $A$ be the intersection $A=X_1\cap X_2$.
If $(X_1,A)$ is $n_1$-connected and
$(X_2,A)$ is $n_2$-connected, with $n_1,n_2\geq 0$,
then for each basepoint $*_V\from V \to A$ with $V\in \V$
the induced map
\[
i_*\from
\pi_r^\V(X_1,A;*_V) \to
\pi_r^\V(X,X_2;*_V)
\]
is an isomorphism if $r<n_1+n_2$  and
is a surjection if $r=n_1+n_2$.
\end{theo}
\begin{remark}
Using the spaces $(X^V,X_1^V,X_2^V,A^V)$ 
and the fact that $A^V$ is a $(\V,\ff^V)$-neighborhood deformation
retract in $X_1^V$ and $X_2^V$ it is possible to 
prove theorem \ref{theo:BMp} also by using the 
adjoints of maps and the homotopy excision theorem of 
Spanier \cite{spanier}.
\end{remark}
\begin{proof}
As in the proof of 16.27 of \cite{gray}
it suffices to consider the case when
$A=X_1\cap X_2$ and $X_1$ and $X_2$ are obtained from $A$ 
by attaching a single cell, namely $X=X_1\cup_A X_2$
with
\[
X_1 = A \cup_\beta(S\times D^m)
\]
\[
X_2= A \cup_\alpha(R\times D^n)
\]
with $n=n_1+1$ and $m=n_2+1$ 
and $R$, $S\in \V$.
In this case we generalize the proof of 13.6 of \cite{gray}
as follows.

Let $p\in e^n$ and $q\in e^m$ and let $R_p = R$
and $S_q = S$ be the fibres of $X$ over $p$ and $q$
respectively.
Then we obtain the commutative diagram
\ldiag{diag:proof}{
\pi_r^\V(X_1,A;*_V) \ar[r]^{i_*} \ar[d]^{\cong} & 
\pi^\V_r(X,X_2;*_V) \ar[d]^{\cong} \\
\pi_r^\V(X\minus R_p,X\minus R_p\minus S_q; *_V) \ar[r]^{\hspace{24pt}j_*} &
\pi_r^\V(X,X\minus S_q; *_V) \\
}%
in which the vertical arrows are isomorphisms. This is readily seen 
by applying 
\ref{lemma:pushout}.
The diagram is the analogue of 13.7 of \cite{gray}.

Next we consider an $\V$-map
\ndiag{%
V\times I^r \ar[r]^{\hspace{-68pt}h'} \ar[d] & 
X = A \cup (R\times E^n) \cup (e^m \cup S) \ar[d] \\
I^r \ar[r]^{h} & \bar X = \bar A\cup e^n \cup e^m
}%
which induces the map $h$ in 13.9 of \cite{gray}.
Using $H_t$ in the proof of 13.9 of \cite{gray}
we obtain a $\V$-homotopy 
\ndiag{%
V\times I^r \ar[d] \ar[r]^{H'_t} & X \ar[d] \\
I^r \ar[r]_{H_t} & \bar X \\ 
}%
by defining for  $(x_1,\dots, x_r)\in I^r$, $v\in V$
\[
H'_t(x_1,\dots,x_3,v) = h'(x_1,\dots,x_{r-1},\lambda,v)
\]
with 
\[
\lambda = 1 - (1-x_r)(1-t\varphi(x_1,\dots,x_{r-1}))
\]
where $\varphi$ is the map in the proof of 13.9 of \cite{gray}.
One readily checks that $(H'_t,H_t)$
is a well-defined $\V$-homotopy.

We can now apply the arguments in the proof of 13.6 of \cite{gray}.
We first show that if $r\leq m+n-2$, $i_*$ is surjective.
Let 
\[
f'\from V\times (I^r,\partial I^r,J^{r-1}) \to  (X,X_2,X_*)
\]
which represents an element $\{f'\}$ in
$\pi^\V_r(X,X_2;*_V)$ as in the proof of 
\ref{lemma:nconnected} above.
Hence $f'$ induces a map
\[
f\from (I^r,\partial I^r, J^{r-1}) \to (\bar X, \bar X_2, *)
\]
with $\bar X = \bar A\cup e^n\cup e^m$.
We obtain for $f$ the homotopy $f_t$ 
as in the proof of 13.6 of \cite{gray} and we define the $\V$-homotopy
\[
f'_t\from V\times (I^r,\partial I^r, J^{r-1}) \to (X,X_2,X_*)
\]
which induces $f_t$ by
\[
f'_t(v,u) =
\begin{cases}
h'_t(v,u) & \mbox{ for $u\in \bar U$, $v\in V$}\\
k'_t(v,u) & \mbox{ for $u\in \bar V$, $v\in V$}\\
f'(v,u) & \mbox{ for $u\in I^r \minus \bar U \minus \bar V$, $v\in V.$}\\
\end{cases}
\]
Here $h'_t$ is defined already in the proof of 
\ref{lemma:nconnected} and $k'_t$ is defined in the same way
as $h'_t$. Hence $f'_t\from f \equiv f'_t$ is a $\V$-homotopy
and $f'_1$ represents $\{f'\}$.
Moreover since $H_1(I^r)\subset \bar X \minus p$
we see that $\{f'\} = \{H'_1\}$ is in the image 
of $j_*$ in \ref{diag:proof}.
This shows that $j_*$, and hence $i_*$, is onto.
In a similar way one follows the argument in the proof 
of 13.6 of \cite{gray}
to show that $i_*$ is injective for $n<n_1+n_2$.
\end{proof}

%\section{The fundamental groupoid of a $\V$-family}
%(It needs a definition without function spaces).
\section{The category $\Pi_\V(X)$ and $\Pi_\V(X)$-modules}
Let $X$ be a $\V$-family and let $x,y\from V\to X$ be $\V$-points
in $X$.
Then we can consider homotopies
$F\from x\equiv y$ with $F\from IV\to X$.
Such a homotopy is termed 
a \emph{$\V$-path} in  $X$.
Homotopy classes of such $\V$-paths relative the boundary
$V\coprod V$ of the cylinder $IV$ 
are termed 
\emph{$\V$-tracks} $\{F\}\from x \tto y$.
Hence a $\V$-track 
is an element
\[
\{F\} \in [IV,X]^{V\coprod V} = [IV,X]^{(x,y)}.
\]
Addition of homotopies yields the composite of tracks $H$, $G$ denoted
by $H\square G$.
Accordingly let $\Pi_\V^V(X)$ be the following groupoid. Objects
are $\V$-maps $V\to X$ and morphisms are $\V$-tracks. For a 
closed inclusion $A\subset X$ 
let 
\begin{equation}
\Pi_\V^V(X,A) \subset \Pi_\V^V(X)
\end{equation}
be the full subgroupoid consisting of objects
$x\from V \to A\subset X$.

The following category $\Pi_\V(X)$ 
plays the role of the fundamental groupoid of a space $X$ in classical
homotopy theory.

\begin{defi}
Let $X$ be a $\V$-family. First we define the 
%\emph{fundamental category} 
category
$\PP_\V(X)$ 
as follows.
Objects are pairs $(V,x)$ where $V$ is an object in $\V$
and 
\[
x\from V\to X
\]
is a $\V$-point in $X$.
A morphism $(V,x) \to (W,y)$ is a 
pair $(a,\alpha)$ where $\alpha\from V \to W$ 
is a morphism in $\V$ and 
$a\from x \tto y\alpha$ 
is a track.
We associate to the morphism $(a,\alpha)$ the diagram
\ndiag{%
V \ar[rr]^{\alpha}\ar[dr]_x &
\ar@{}[d]|{\displaystyle\stackrel{a}{\tto}}
& W \ar[dl]^{y} \\
& X \\
}%
The composition of such morphisms
according to the diagram
\ndiag{
U \ar[d]^z \ar[r]^{\beta}
\ar@{}[dr]|{\displaystyle\stackrel{b}{\tto}}
& V \ar[d]^x \ar[r]^{\alpha} 
\ar@{}[dr]|{\displaystyle\stackrel{a}{\tto}}
& W \ar[d]^y \\
X \ar@{=}[r] & X \ar@{=}[r] & X \\
}%
is defined by
\[
(a,\alpha)(b,\beta) = ((a\beta) \square b, \alpha \beta).
\]
\end{defi}

Two morphisms $(a,\alpha)$ and $(a',\alpha')\from (V,x) \to (W,y)$
are \emph{equivalent} if there is a track
\[
h\from \alpha \tto \alpha' \mbox{ \ \ in $W^V$}
\]
such that $(yh)\square a = a'$.
This is a natural equivalence relation on $\PP_\V(X)$ 
so that we obtain the quotient category
\begin{equation}\label{eq:fundcat}
\Pi_\V(X) = \PP_\V(X)/_\sim
\end{equation}
which is termed the (discrete) \emph{fundamental category} of $X$.

For a closed inclusion $A\subset X$ in $\Vtop$ let 
\begin{equation}
\label{eq:inclusion}
\Pi_\V(X,A) \subset \Pi_\V(X)
\end{equation}
be the full subcategory consisting of all objects
$(V,x)$ with $x\from V\to A\subset X$.

The cellular approximation theorem shows for
a $0$-connected relative
$\V$-complex $(X,A)$ 
that the inclusion $X_2\subset X$ induces an isomorphism
\begin{equation}
\label{eq:isomorphism}
\Pi_\V(X_2,A) = \Pi_\V(X,A).
\end{equation}

\begin{remark}
\label{rem:5.2}
For a space $X$ in $\top$ let $\Pi(X)$ be the fundamental groupoid
of $X$. Then we get the functor 
\begin{equation*}
\Pi (X^\circ)\from \V^\op \to \grd
\end{equation*}
which carries $V$ to the fundamental groupoid $\Pi(X^V)$ of the function
space $X^V$. 
We have $\Pi(X^V) = \Pi_\V^V(X)$.
Now the category $\PP_\V(X)$ above coincides with the 
``integration category''
$\int_\V \Pi(X^\circ)$,
see \cite{mosv93} and compare with \cite{lueck}.
\end{remark}

A left (resp. right) \emph{$\Pi_\V(X)$-module}
is a covariant (resp. contravariant) functor
\[
M\from \Pi_\V(X) \to \Ab
\]
where $\Ab$ is the category of abelian groups.
Let $\Mod(\Pi_\V(X)^\op)$ be the category of right modules. Morphisms
are natural transformations.

\begin{propo}
\label{propo:1}
For $n\geq 2$ the homotopy group
yields a right module
\[
\pi_n^\V(X)\from \Pi_\V(X)^\op \to \Ab.
\]
Moreover for $n\geq 3$ relative homotopy groups of a pair $(Y,X)$
yield the right module
\[
\pi_n^\V(Y,X)\from \Pi_\V(X)^\op \to \Ab.
\]
\end{propo}
\begin{proof}
An object $(V,x)$ of $\Pi_\V(X)$ yields 
a basepoint $x\from V \to  X$, hence the abelian
group $\pi_n^\V(X;x)$ is defined for $n\geq 2$.
If $(a,\alpha)\from (V,x) \to (W,y)$ is a morphism
in $\PP_\V(X)$ 
and $n\geq 2$
then
there is a canonical
induced homomorphism of abelian groups defined by 
the composition
\ldiag{diag:hom}{
\pi_n(X^W;y) \ar@{=}[d] &
\pi_n(X^V;y\alpha) \ar@{=}[d] &
\pi_n(X^V;x) \ar@{=}[d]\\
\pi_n^\V(X;y) \ar[r]^{\pi_n(\alpha^*)} &
\pi_n^\V(X;y\alpha) \ar[r]^{a^\#} &
\pi_n^\V(X;x)\\
}%
where
\[
\alpha^*(y)=y\alpha \from \xymatrix@1{V \ar[r]^\alpha & W \ar[r]^y & X}
\]
and $a^\#$ is the isomorphism of homotopy
groups induced by the track $a$. 
It is not difficult to check that 
the structure
of $\PP_\V(X)$-module induces a structure of $\Pi_\V(X)$-module,
since the homomorphism  in \ref{diag:hom} does not depend on
the representative $(a,\alpha)$ in $\PP_\V(X)$ of a morphism
in $\Pi_\V(X)$.
If $n\geq 3$ the same argument can be applied to the relative group
\(
\pi_n^\V(Y,X;x).
\)
\end{proof}

Let $Z$ be a $\V$-set and let $\alpha\from Z \to X$ 
be a $\V$-map.
Hence $(Z,\alpha)$ is given by a family of $\V$-maps
$\alpha|V\from V \to X$ with $V\in Z$.
We define the \emph{free right $\Pi_\V(X)$-module} in
$\Mod(\Pi_\V(X)^\op)$
\begin{equation}
\label{eq:free}
M[Z,\alpha]\from 
\Pi_\V(X)^\op \to \Ab
\end{equation}
by the direct sum
\[
M[Z,\alpha] = 
\bigoplus_{V\in Z}
\ze\hom(-,(V,\alpha|V)).
\]
Here $\ze\hom(-,(V,\alpha|V))$ denotes the module which
carries $(W,y)$ to the free abelian
group generated by the hom-set
$\hom((W,y),(V,\alpha|V))$ of 
all morphisms $(W,y) \to (V,\alpha|V)$ in $\Pi_\V(X)$.
A morphism
\begin{equation}
\varphi\from M[Z,\alpha] \to M
\end{equation}
in $\Mod(\Pi_\V(X)^\op)$
is uniquely determined by a family of elements
\[
\varphi_V \in M(V,\alpha|V) 
\mbox{\ \ for $V\in Z$}
\]
such that 
$\varphi(1_{(V,\alpha|V)}) = \varphi_V$.
Here $1_{(V,\alpha|V)}$ denotes the identity of the object
$(V,\alpha|V)$ in $\Pi_\V(X)$.

We now consider for $(Z,\alpha)$ above the pushout
\ndiag{
Z\times * \ar@{=}[r] \ar[d] & Z \ar[r]^{\alpha} & X \ar[d] \\
Z\times S^n \ar@{=}[r] & S^n_Z  \ar[r] & X\cup_\alpha S^n_Z \\
}%

The projection $Z\times S^n \to Z$ induces the retraction
\[
r_X\from X\cup_\alpha S^n_Z \to X
\]
which is the identity on $X$. For $n\geq 1$ let 
\[
\pi^\V_n(X\cup_\alpha S^n_Z)_X = 
\ker \left(
(r_X)_*\from \pi_n^\V(X\cup_\alpha S^n_Z) \to \pi_n^\V(X)
\right).
\]
We define the \emph{partial suspension}
\begin{equation}
E\from 
\pi_n^\V(X\cup_\alpha S^n_Z)_X \to 
\pi_{n+1}^\V(X\cup_\alpha S^{n+1}_Z)_X
\end{equation}
by the composition $E=j^{-1}\pi_* \partial^{-1}$
\ndiag{%
\pi_{n+1}^\V(X\cup_\alpha D_Z^{n+1},X \cup_\alpha S^n_Z)
\ar[r]^{\hspace{24pt}\partial}_{\hspace{24pt}\cong} 
\ar[d]_{\pi_*}
& 
\pi_n^\V(X\cup_\alpha S^n_Z)_X \\
\pi_{n+1}^\V(X\cup_\alpha S^{n+1}_Z,X) & 
\pi_{n+1}^\V(X\cup_\alpha S^{n+1}_Z)_X \ar[l]^{j}_{\cong} \\
}%
Here $\pi$ is induced by $D^{n+1}/S^{n} = S^{n+1}$ and the isomorphisms
$\partial$ and $j$ are induced by the homotopy exact sequences
of pairs. The Blakers-Massey theorem implies:

\begin{propo}
$E$ is an isomorphism for $n\geq 2$ and is surjective for $n=1$.
\end{propo}

\begin{propo}
\label{propo:2}
There is an isomorphism of modules $(n\geq 2)$
\[
\varphi\from M[Z,\alpha] \cong \pi_n^\V(X\cup_\alpha S^n_Z)_X
\]
which is given by the family of maps 
\[
\varphi_V\from S^n_V\subset 
S^n_Z \to X\cup_\alpha S^n_Z
\]
with $V\in Z$. Moreover $\varphi$ is compatible with $E$.
\end{propo}
\begin{proof}
Without loss of generality, by additivity, we can assume that 
$Z=W$; moreover, by considering
the mapping cylinder of $\alpha\from W \to X$ 
in $\Vtop$ 
we can assume that $\alpha\from W \to X$ is a closed inclusion
(actually, a $\V$-cofibration). 
Let $(V,x)$ be an arbitrary object in $\Pi_\V(X)$. 
We want to show that $\varphi$ induces an isomorphism
of groups
\begin{equation}\label{eq:varphi}
\varphi\from  M[W,\alpha](V,x) = \ze\hom( (V,x), (W,\alpha)) 
\to
\pi^\V_n(X\cup_\alpha S^n_W;x)_X.
\end{equation}

Let $X_x^V$ denote the path-component in $X^V$ containing
$x\in X^V$ and $W_x^V = X^V_x \cap \alpha^V(W^V)$;
let $\alpha_x\from W^V_x\subset X^V_x$ be the inclusion.
Moreover, let $\widetilde{X_x^V}=E_X$ be the universal covering
space of $X_x^V$.
We can assume that any space  that we 
are considering in this proof has a universal cover by
taking a CW-approximation weakly equivalent to it.
The elements of $\widetilde{X_x^V}$ correspond bijectively
with tracks $x \tto \xi$, where $\xi\from V \to X$  
is any $\V$-map. The covering projection
in this case is the evaluation of the track at $1$.
We can define $E_W$ 
by the following pull-back diagram.
\ldiag{diag:pull}{%
E_W \ar@{->>}[d] 
\ar@{ >->}[r]^{\hspace{-12pt}\widetilde{\alpha_x}} & \widetilde{X^V_x} = E_X\ar@{->>}[d] \\
W^V_x \ar@{ >->}[r]^{\alpha_x} 
& X^V_x  \\
}%

Then $E_X\cup_{\widetilde{\alpha_x}} ( E_W\times S^n)$
is the universal covering of $X^V_x\cup_{\alpha_x} (W_x^V\times S^n)$.

Consider the morphisms $(b,\beta)$ 
in $\PP_\V(X)( (V,x), (W,\alpha))$.
Then $\beta$ and $b$ are elements in
%The projections onto the first and the second factor of $(b,\beta)$ 
%yield two maps to 
$\widetilde{X_x^V}$ and $W^V_x$ 
respectively, 
compatible with the pull-back diagram  
\ref{diag:pull}. 
Hence there is a map $\lambda$ 
as in the following diagram.
\[
\xymatrix@C+36pt{
\PP_\V(X)( (V,x), (W,\alpha))  \ar[r]^\lambda \ar@{->>}[d] &
E_W \ar@{ ->>}[d] \\
\Pi_\V(X)( (V,x), (W,\alpha)) \ar@{.>}[r]^{\bar\lambda}_{\cong} & 
\pi_0(E_W).
}
\]
Here two morphisms $(b,\beta)$ and $(b',\beta')$ in
$\PP_\V(X)( (V,x), (W,\alpha))$
have images under $\lambda$ in the same path-component
if and only if there is $h\from I \to E_W$
such that $h(0) = b$, $h(1) = b'$.
But this happens if and only if $(b,\beta)$ and 
$(b',\beta')$ belong to the same
equivalence class in $\Pi_\V(X)$, and thus 
there is an induced bijection 
$\bar\lambda$ as in the diagram.

Now, using the exact homotopy sequence of the pair,
the property of universal covering spaces, 
the homotopy excision, and the isomorphism induced by $\bar\lambda$
we get the equations:
\begin{equation*}
\begin{split}
\pi^\V_n(X\cup_\alpha S^n_W;x)_X &=
\pi^\V_n(X\cup_\alpha S^n_W,X;x) \\
&= 
\pi_n(X^V\cup_{\alpha^V} S^n_{W^V},X^V;x) \\
&=
\pi_n(X_x^V\cup_{\alpha_x^V} S^n_{W_x^V},X^V_x;x) \\
&=
\pi_n(E_X\cup_{\widetilde{\alpha_x}} (E_W\times S^n),
E_X;x) \\
&=
\pi_n(E_X\cup_{\widetilde{\alpha_x}} (E_W\times S^n),
E_X) \\
&\cong
\bigoplus_{C\subset E_W}
\pi_n(S^n\times C,C) \\
&=
\ze\pi_0(E_W) \\
&=
\ze\Pi_\V(X)( (V,x), (W,\alpha)) 
\end{split}
\end{equation*}
where
the sum $\bigoplus$ ranges over the set of path-components
$C\subset E_W$. 
We can apply the homotopy excision since $E_X$ is 
simply connected and $W=W\times *$ is a $\V$-neighborhood deformation
retract in $S^n_W$.
Therefore we have the isomorphism 
\begin{equation*}
 M[W,\alpha](V,x) = \ze\hom( (V,x), (W,\alpha)) 
\cong
\pi^\V_n(X\cup_\alpha S^n_W;x)_X.
\end{equation*}
It is not difficult, by chasing along the chain of equalities above,
to see that this is the same homomorphism induced by $\varphi$ 
as in  \ref{eq:varphi} and that it commutes with the partial
suspension $E$
as claimed. Since $\varphi$ is a morphism
of $\Mod(\Pi_\V(X)^\op)$-modules, the proof is complete.
\end{proof}

\section{Homology and cohomology}

Let $(X,A)$ be a relative $\V$-complex. We say that $(X,A)$ 
is \emph{normalized} if the attaching maps
\begin{equation}
f_n\from S^{n-1}_{Z_n} = Z_n \times S^{n-1} \to X_{n-1}
\end{equation}
of $n$ cells in $X$ ($n\geq 1$) carry the basepoint $*$ of $S^{n-1}$
to the $0$-skeleton, that is $f_n(Z_n\times *) \subset X_0$
or equivalently if the attaching maps are cellular.
Moreover $(X,A)$ is \emph{reduced}
if $X_0=A$ so that $X$ is obtained from $A$ 
by
attaching $\V$-cells of dimension $\geq 1$.

\begin{lemma}\label{lemma:rednorm}
Let $(X,A)$ be a relative $\V$-complex which is $0$-connected.
Then there exists a relative $\V$-complex $(Y,A)$
which is reduced and normalized together with a $\V$-homotopy
equivalence $Y\to X$ under $A$.
\end{lemma}
\begin{proof}
First we obtain a reduced $\V$-complex 
$(Y',A)$ by \ref{lemma:approximation}.
Using the cellular approximation theorem \ref{theo:cat}
the attaching maps in $Y'$ are $\V$-homotopic 
to normalized attaching maps.
This can be used to construct inductively $Y$.
\end{proof}

If $(X,A)$ is reduced and normalized we obtain for $n\geq 1$
the commutative diagram
\diag{%
Z_n\times S^{n-1} \ar@{=}[r] & S^{n-1}_{Z_n} \ar[r]^{f_n} & X_{n-1} \\
Z_n\times * \ar@{=}[r] \ar[u] & Z_n \ar[r]^{\alpha_n} & A \ar[u] \\
}%
where $\alpha_n$ is the \emph{basepoint map}
associated to the attaching map $f_n$.
We identify the closed cell $D^n$ with the 
(reduced) cone over the pointed space $S^{n-1}$ so that we obtain the pinch 
map under $S^{n-1}$ ($n\geq 1$)
\[
\xymatrix@1{D^n \ar[r]^{\hspace{-24pt}\bar \mu} & D^n \vee S^{n}.} 
\]
This map induces the coaction map
\diag{%
X_n \ar[r] ^{\mu} \ar@{=}[d] & X_n\cup_{\alpha_n} S^n_{Z_n} \ar@{=}[d] \\
X_{n-1}\cup_{f_n}(Z_n\times D^n) \ar[r] & X_{n-1}\cup_{f_n}(Z_n\times(D^n\vee S^n)) \\
}%
where the bottom row is $1_{X_{n-1}} \cup(Z_n\times \bar \mu)$.
On the other hand we have for the $1$ dimensional disk $D^1$ with 
boundary points $\partial D^1 = \{\partial_0,\partial_1\}$ 
and basepoint $\partial_0=*$ the map
\ndiag{%
(S^1,\partial_0) \ar[r]^{\hspace{-24pt}\bar \mu_0} 
& (D^1\cup_{\partial_1} S^1, \partial_0) \\
}%
which is defined by composing obvious tracks $\partial_0 \to \partial_1$ 
in $D^1$, 
$*\tto *$ in $S^1$ and 
$\partial_1 \tto \partial_0$ 
in $D^1$.
If $X$ is a normalized $\V$-complex we have the $\V$-set $X_0=Z_0$ 
and the inclusion $\alpha_0\from Z_0 \subset X$.
Then $\bar \mu_0$ induces the map
\ldiag{diag:6.5}{
S^1_{Z_1} \ar[r]^{\mu_0} & X\cup_{\alpha_0} S^1_{Z_0} \\
Z_1\times S^1 \ar@{=}[u] \ar[r] & 
Z_1 \times (D^1\cup _{\partial_1} S^1) \ar[u]\\
}%
Here the bottom row is $Z_1 \times \bar \mu_0$ and the right hand side is
$\bar f_1 \cup (S^1\times \alpha_1^1)$ 
where $\bar f_1$ is the characteristic map of $1$-cells
and $\alpha_1^1 = f_1|(Z_1\times \partial_1)$.
The map $\mu_0$ restricted to $Z_1\times *$ is given by
$\alpha_1 = f_1|(Z_1\times \partial_0)$.

\begin{defi}
Let $(X,A)$ be a $\V$-complex which is reduced and normalized.
Then we define the \emph{chain complex} $C_*(X,A)$ 
in $\Mod(\Pi_\V(X,A)^\op)$ as follows.
For $n\leq 0$ let $C_n(X,A) = 0 $ and for 
$n\geq 1$ let 
\[
C_n(X,A) = M[Z_n,\alpha_n]
\]
be defined as in \ref{eq:free} by the $\V$-set $Z_n$ 
of $n$-cells in $(X,A)$ and the basepoint map
$\alpha_n\from Z_n \to A\subset X$.
The differential 
\[
d\from C_{n+1}(X,A) \to C_n(X,A)
\]
is defined for $V\in Z_{n+1}$ by elements $d_V$ obtained as follows.
Let $\partial_V$ be the composite ($n\geq 1$)
\[
\partial_V\from 
\xymatrix@1@R+24pt{ 
S^n_V \subset S^{n}_{Z_{n+1}} \ar[r]^{\hspace{12pt}f_{n+1}} & 
X_n \ar[r]^{\hspace{-58pt}\mu} & 
X_n \cup_{\alpha_n} S^n_{Z_n} \subset 
X \cup_{\alpha_n} S^n_{Z_n} 
}
\]
which for $n\geq 2$ represents the element 
\[
d_V\in \pi^\V_n(X\cup_{\alpha_n} S^n_{Z_n}; \alpha_{n+1}|V)_X 
=
M[Z_n,\alpha_n](V,\alpha_{n+1}|V).
\]
Compare \ref{rem:5.2} and \ref{eq:free}.
For $n=1$ we have the partial suspension
\[
E\from \pi_1^\V(X\cup_{\alpha_1} S^1_{Z_1}; \alpha_2|V)_X \to 
\pi_2^\V(X\cup_{\alpha_1} S^2_{Z_1}; \alpha_2|V)_X
\]
and we set $d_V = E(\partial_V)$.
\end{defi}

If $X$ is a normalized $\V$-complex we define the chain complex $C_* X $ 
in $\Mod( \Pi_\V(X,X_0)^\op)$ 
by
\begin{eqnarray*}
C_n(X) &=& C_n(X,X_0) \mbox{ \ \ for $n\geq 1$ and } 
\\
C_0(X) &=& M[Z_0,\alpha_0] 
\end{eqnarray*}
where $Z_0 = X_0$ is the $0$-skeleton of $X$ which is an $\V$-set
and $\alpha_0\from Z_0 \subset X$ is the inclusion.
Moreover the differential
\[
d\from C_{n+1}(X) \to C_n(X) 
\]
is defined as above for $n>1$. For $n=1$ the differential
\[
d\from C_1(X) \to C_0(X)
\]
is defined by $d_V = E(\mu_0|S^1_V)$ for $V\in Z_1$ 
where $\mu_0$   is the map in \ref{diag:6.5}
and $E$ is the partial suspension. 
We may consider $C_*(X,A)$ and $C_*(X)$ also as chain complexes
in $\Mod(\Pi_\V(X)^\op)$.

\begin{propo}
$C_*(X,A)$ and $C_*(X)$ are well defined chain complexes.
\end{propo}
\begin{proof}
Since $(\Vtop,\TT)$ is homological
\ref{rem:pro},  the \emph{chain functor} $C_*$ 
is defined as in V.2.4, page 255 of \cite{baues},
and it is easy to see that it coincides with 
the chain functor $C_*$ defined above.
% It is only left to show that the coefficient  \todo
\end{proof}

\begin{defi}
\label{defi:homo}
Let $(X,A)$ be a $0$-connected relative $\V$-complex. 
By lemma \ref{lemma:rednorm}
we can assume that $(X,A)$ is reduced and normalized. 
Given a \emph{right} $\Pi_\V(X)$-module $M$ we define the cohomology
of $(X,A)$ with coefficients in $M$ by
\[
H^*(X,A;M) = H^*\left( \hom(C_*(X,A),M ) \right),
\]
where $\hom$ is defined in the abelian category
of right $\Pi_\V(X)$-modules.

In a similar way, if $N$ is a left $\Pi_\V(X)$-module
we can define the \emph{homology} of $(X,A)$ with coefficients
in $N$ by
\[
H_*(X,A;N) = H_* ( C_*(X,A)\otimes N ),
\]
where $\otimes$ stands for the tensor product
of a right and a left $\Pi_\V(X)$-module.
Here $H^*(X,A;M)$ and $H_*(X,A;M)$ are abelian groups.

Finally, the \emph{total homology}
of $(X,A)$ is defined by
\[
H_*(X,A) = H_* (C_*(X,A) ).
\]
The total homology is the homology of a chain complex of right
$\Pi_\V(X)$-modules, so that $H_*(X,A)$ is again
a $\Pi_\V(X)$-module.
\end{defi}

\begin{remark}
This is a generalization of the 
Bredon--Illman cohomology \cite{bredon} for $G$-spaces.
See Br\"ocker for the singular Bredon homology \cite{brocker}.
See \cite{moller,mosv95,mosv93} for singular Bredon cohomology
of discrete diagrams.
\end{remark}

\section{The Whitehead theorem}
\label{sec:whitehead}
\begin{theo}
\label{theo:whitehead2}
Let $A$ be a $\V$-family and let $f\from (X,A) \to (Y,A)$
be a cellular map between normalized reduced
relative $\V$-complexes. Then $f$ is a $\V$-homotopy equivalence under $A$
if and only if $f$ induces an isomorphism 
\[
\varphi=\Pi(f)\from \Pi_\V^V(X,A) \to \Pi_\V^V(Y,A)
\]
of groupoids for all $V\in \V$
and any one of the following cases is true:
\begin{enumerate}
\item
$f$ induces a homotopy equivalence of $\Pi_\V(X,A)$-chain complexes
\[
f_*\from C_*(X,A) \to C_*(Y,A),
\]
where we use the isomorphism $\varphi$ 
to identify $\Pi_\V(X,A)$ and $\Pi_\V(Y,A)$.
\item
The induced map between homology
\[
f_*\from H_*(X,A) \to \varphi^*H_*(Y,A)
\]
is an isomorphism.
\item
For every $\Pi_\V(Y,A)$-module $M$ the induced map 
\[
f^*\from H^*( Y,A; M) \to H^*(X,A;\varphi^* M)
\]
is an isomorphism.
\end{enumerate}
\end{theo}
\begin{proof}
This is a consequence of the general homological
Whitehead theorem 
V.7.1 of \cite{baues}.
It is only needed  to show that any reduced and normalized $\V$-complex
is a $\mathbf{T}$-complex in the sense of IV.2.2 of \cite{baues},
where $\mathbf{T}$ the theory of $\V$-graphs (see definition \ref{defi:T}). 
This follows from lemma \ref{lemma:Tcomplex} in the appendix.
\end{proof}

\begin{proof}[Proof of theorem \ref{theo:whitehead}]
By \ref{lemma:nconnected} the $\TT$-complexes in $\Vtop$
are $\TT$-good in the sense of IV.3.7 \cite{baues}.
Hence by \ref{lemma:rednorm} and IV.3.11 of \cite{baues}
we know that the Whitehead theorem  holds for $0$-connected relative
$\V$-complexes $(X',A)$ and $(Y',A)$.
This can be used to prove the general case in \ref{theo:whitehead}
as follows. Let $f\from (X,A) \to (Y,A)$ be a Whitehead equivalence
as in \ref{theo:whitehead}. 
We may assume that $f$ is cellular. Hence $f$ yields a restriction
$\alpha\from X_0\to Y_0 \subset Y$  which is surjective in $\pi_0^\V$.
Let $M_\alpha = IX_0\cup_\alpha Y$ 
be the mapping cylinder (in $\Vtop$) of $\alpha$, i.e. 
the push-out of $i_0\from X_0 \to IX_0$ and $\alpha$ (see \ref{coro:factorization}).
Then $i_1\from X_0 \to IX_0$ 
yields a pair $(M_\alpha,X_0)$ which is $0$-connected.
Moreover $Y \to M_\alpha $ is a $\V$-homotopy equivalence
and the composite $X\to Y \to M_\alpha$
is by the homotopy extension property of $X_0\subset X$ 
$\V$-homotopic to a map $g\from X \to M_\alpha$ 
which is the identity on $X_0$. Hence 
$g\from (X,X_0) \to (M\alpha,Y_0)$ 
is a Whitehead equivalence between $0$-connected pairs and hence
a $\V$-homotopy equivalence.
\end{proof}

\section{The Whitehead sequence and obstruction theory}
Let $(X,A)$ a reduced normalized $\V$-complex and 
$\Pi_\V(X,A)$ the (restricted) fundamental category.
Then the homotopy groups $\pi_n^\V(X)$ and the homology groups 
$H_n(X,A)$ 
are 
$\Pi_\V(X,A)^\op$-modules and it is possible to define 
a Hurewicz homomorphism 
\[
h\from \pi_n^\V(X) \to H_n(X,A)
\]
as follows.

Let $x\from V\to X$ be a basepoint and  $\alpha\from S^n_V \to X$ 
a representative of the element $[\alpha] \in \pi^\V_n(X;x)$.
By \ref{theo:cat} we can assume that $\alpha S^n_V \subset X_n$.
Consider the difference 
$\mu\alpha - i\alpha$
of the two composites
\ndiag{%
S^n_V \ar[r]^{\alpha} \ar[d]^\alpha & X_n \ar[d]^{\mu} \\
X_n \ar[r] ^{\hspace{-36pt}i} & X_n\cup_{\alpha_n} S^n_{Z_n} \ar@{ >->}[r]^{j} & 
X\cup_{\alpha_n} S^n_{Z_n} \\
}%
where $\mu$ is the coaction map.
Then  $j(\mu\alpha - i \alpha)$
is a cycle in $C_n(X,A)$ that represents an element 
\[
h(\alpha) = \left\{j(\mu\alpha - i\alpha)\right\} \in H_n(X,A).
\]
It is not difficult to see that $h(\alpha)$ does not depend
upon the $\V$-homotopy class of $\alpha$ and that $h$ 
is a $\Pi_\V(X,A)$-homomorphism, termed the \emph{Hurewicz}
homomorphism.

The groups $\Gamma_n(X,A)$, if $n\geq 3$, are defined as follows.
If $x\from V\to A$ is an object of $\Pi_\V(X,A)$ then 
\[
\Gamma_{n}(X,A)(x) = 
\image \left[ 
\pi_n^\V(X_{n-1};x) \to \pi_n^\V(X_{n};x)
\right]
\]
For $n=1,2$ the definition is more complicate 
(see  section V.5, page 262, of \cite{baues}).

\begin{theo}
\label{theo:hurewicz}
Let $(X,A)$ be a reduced and normalized $\V$-complex 
Then the following sequence of 
$\Pi_\V(X,A)^\op$-modules is exact
\begin{equation*}
\begin{split}
\xymatrix@1{%
\ar[r] & 
\Gamma_n(X,A) \ar[r] & 
\pi_n^\V(X) \ar[r]^h & 
H_n(X,A) \ar[r] & 
\Gamma_{n-1}(X,A) \ar[r] & 
\dots 
}
\\
\xymatrix@1{%
\ar[r] & 
\Gamma_2(X,A) \ar[r] & 
\pi^\V_2(X) \ar[r]^h &
H_2(X,A) \ar@{->>}[r] & 
\Gamma_1(X,A) \ar[r] & 
0
}
\end{split}
\end{equation*}
Furthermore, the sequence is natural in $(X,A)$ in the category 
$\Vtop^A$ of $\V$-families under $A$. 
\end{theo}
\begin{proof}
It is a direct consequence of theorem V.5.4, page 264 \cite{baues}.
\end{proof}

%\begin{comment}
%\begin{remark}
%In the framework of equivariant stable homotopy theory
%\cite{lmsm},
%Lewis defined a Hurewicz homomorphism in \cite{lewis}
%for equivariant ($RO(G)$-graded) homotopy groups and 
%homology. 
%Since for a (compact) Lie group $G$ a $G$-space is 
%exactly a $\V$-family, where $\V$ is the orbit category
%of $G$, one might try to compare the homotopy 
%and the homology with Burnside ring coefficients 
%in \cite{lmsm} 
%with the corresponding $\V$-homotopy and homology of this paper.
%%But actually the equivariant homotopy theory of \cite{lmsm} 
%%has such a different structure, that it is probably 
%%unreasonable to try to relate the $RO(G)$-graded theory
%%with our $\V$-homotopy theory
%%(which under some extent 
%%is closer to the equivariant algebraic $K$-theory of \cite{lueck}).
%\end{remark}
%\end{comment}

The cohomology defined above is suitable for obstruction theory
as follows.

% Obstruction theory
\begin{theo}
Let $(X,A)$ a normalized reduced relative $\V$-complex and $f\from A \to Y$
a $\V$-map that admits an extension $g\from X_n\to Y$,
for $n\geq 2$. Then its restriction $g|X_{n-1}$ has an extension
$g'\from X_{n+1} \to Y$
to the $(n+1)$-skeleton 
if and only if the obstruction 
\[
\obst(g|X_{n-1}) \in H^{n+1}(X,A; g^*\pi^\V_nY)
\]
vanishes, where $g_*\from \Pi_\V(X,A) \to \Pi(Y)$ 
induces the functor $g^*\from \Mod(\Pi_\V(Y)) \to \Mod(\Pi_\V(X,A))$.
\end{theo}
\begin{proof}
The obstruction $\obst(g|X_{n-1})$
is defined as  in chapter V of \cite{baues}.
See theorem V.4.4, page 262.
\end{proof}

\section{Wall finiteness obstruction}
% domination
A $\V$-complex $X$ is \emph{finite} if it has only finitely
many cells or equivalently if $\bar X$ is a finite CW-complex.
A \emph{domination} of $Y$ in $\Vtop$ is given
by $\V$-maps
\[
\xymatrix@1{ Y \ar[r]^f & 
X \ar[r]^g &
Y \\}
\]
and a $\V$-homotopy $H\from gf \sim 1_Y$.
The domination has dimension $\leq n$ if the dimension 
of $X$ is $\leq n$ and the domination is finite if $X$ is a finite
$\V$-complex.
% projective class group
% theo
% lemma
For a ringoid $R$ let $K_0(R)$ denote the reduced projective class
group of $R$.
If $\CC$ is a category let $\ze\CC$ be the ringoid associated to $\CC$
for which 
\[
\hom_{\ze\CC}(X,Y) = \ze[ \hom_\CC(X,Y)]
\]
is 
the free abelian group generated by $\hom_\CC(X,Y)$.

\begin{theo}
\label{theo:finobst}
Let $Y$ be a $\V$-complex which admits a finite domination
in $\Vtop$. Then the finiteness obstruction 
\[
[Y] = [C_*(Y)] \in K_0(\ze\Pi_\V(Y,Y_0))
\]
is defined, 
%where $\ze\Pi_\V(Y,Y_0)$ 
%is the enveloping ringoid,
with the property that 
the obstruction $[Y]=0$  is trivial if and only if 
$Y$ is $\V$-homotopy equivalent to a finite $\V$-complex.
Moreover, if the domination of $Y$ has dimension $\leq n$
and $[Y] = 0$ then $Y$ is $\V$-homotopy equivalent to a finite $\V$-complex
of dimension $\leq \mathrm{Max}(3,n)$.
\end{theo}

\begin{remark}
A result like \ref{theo:finobst}
was proved by different methods in the case of 
$G$-spaces in theorem 14.6 of L\"uck \cite{lueck},
without the dimension estimate.
\end{remark}

\begin{lemma}
\label{lemma:obst}
Let $\tilde X$ be a finite domination of $Y$. Then there exists a finite
domination $X$ of $Y$ for which $f\from Y \to X$
induces an isomorphism $\pi_0^\V(Y;*_V) \cong \pi_0^\V(X;f*_V)$
for every basepoint $*_V$.
\end{lemma}
\begin{proof}
We may assume that $\tilde f\from Y \to \tilde X$ 
and $\tilde g\from \tilde X \to Y$
are cellular maps which yield the composite
$\lambda = i \tilde f_0 \tilde g_0\from \tilde X_0 \to \tilde X_0
\subset \tilde X$.
Let $X$ be obtained by attaching $1$-cells to $\tilde X$ 
as in the push-out diagram
\ndiag{
\tilde X_0 \times \partial I \ar[r] \ar[d]^{(i,\lambda)} 
& \tilde X_0 \times I \ar[d] \\
\tilde X \ar[r]^{j} & X \\
}%
where $i\from \tilde X_0 \to \tilde X$ is the inclusion.
Since $\tilde g i \sim \tilde g\lambda$, there exists
an extension $g\from X \to Y$ 
of $\tilde g$ with $gj = \tilde g$.
Hence for $f = j \tilde f$
we obtain the domination
\[
\xymatrix@1{ Y \ar[r]^f & X \ar[r]^g & Y \\}
\]
with $gf = \tilde g \tilde f \sim 1_Y$.
This shows that $f_*\from \pi_0^\V(Y;*_V) \to \pi_0^\V(X;f*_V)$
is surjective.
In fact, $f$ is also surjective since an element of $ \pi_0^\V(X;f*_V)$
can be represented by a $\V$-map $\xi\from V \to X_0=\tilde X_0$
via the cellular approximation theorem \ref{theo:cat}.
For the inclusion $j_0\from X_0 \subset X$ we get the following homotopy
\[
j_0 = ji\xi \sim j\lambda \xi = jif_0g_0\xi = 
j_0f_0(\eta) = f_*(\eta),
\]
where $\eta = g_0\xi$.
\end{proof}

\begin{proof}[Proof of \ref{theo:finobst}]
Let $X$ be a finite domination of $Y$ as in 
\ref{lemma:obst}.
We can choose $tg$ to be cellular and we obtain by restriction of $g$
the map $\alpha = ig_0 \from X_0 \to Y_0\subset Y$.
We use $\alpha$  to construct 
the mapping cylinder 
$Y'=M_\alpha = Y \cup_\alpha IX_0$,
so that $(M_\alpha,X_0)$ 
is a $0$-connected relative $\V$-complex
and $M_\alpha \equiv Y$ in $\Vtop$.
%Here we  use the assumption on $\pi_0$ of \ref{}.\todo
The domination $X$ of $Y$ now yields 
a domination $X'$ of $M_\alpha$
by the composites
\[
\xymatrix@1{%
Y' \equiv Y \ar[r]^{f} & 
X \equiv X' \equiv X \ar[r]^{g}
&
Y \equiv Y'
\\
}%
\]
where the composites $f'\from Y' \to X'$ 
and $g'\from X' \to Y'$ can be chosen
to be maps under $X_0$ 
and where the $\V$-homotopy $g'f'\sim 1_{Y'}$
can be chosen to be a $\V$-homotopy relative $X_0$.
Moreover we can assume by the approximation lemma \ref{lemma:rednorm}
that $(Y',X_0)$ 
and $(X',X_0)$ 
are replaced by normalized and reduced relative $\V$-complexes
$(Y'',X_0)$ and $(X'',X_0)$.
Hence we can apply theorem VII.2.5 of \cite{baues} which shows 
that a finiteness obstruction
\[
[Y''] = [C_*(Y'',X_0)] \in  K_0(\ze\Pi_\V(X,A))
\]
is defined and has the property that $[Y''] = 0$
if and only if there exists a finite relative
$\V$-complex $(Z,X_0)$
which is reduced and normalized 
with $\dim(Z) \leq \mathrm{Max}(3,n)$,
where $n= \dim(X'')$, and such that 
$Z$ is $\V$-homotopy equivalent to $Y''$ relative $X_0$.
Since we have homotopy equivalences of chain complexes
$C_*(Y) \equiv C_*(Y') \equiv C_*(Y',X_0) \equiv C_*(Y'',X_0)$,
and since by the finiteness of $X_0$ 
we have 
$[C_*(Y')] = [C_*(Y',X_0)]$,
we see that $[C_*(Y)] = [C_*(Y'',X_0)]$.
This yields the results in \ref{theo:finobst}.
\end{proof}

\section{Whitehead torsion}
For $n\geq 0$   let $\R_+^{n+1}$  
and $\R_-^{n+1}$ be the subspaces 
of $\R^{n+1}$ defined by the elements $(x_0,\dots,x_n)\in \R^{n+1}$ 
with $x_0\geq 0$ and $x_0\leq 0$ respectively.
Let $D^{n+1}$ denote the unit sphere in $\R^{n+1}$ and 
$S^{n}$ its boundary.
A \emph{ball pair} is a tuple
$(\square^{n+1},S^{n},P^n,Q^n)$ homeomorphic to the Euclidean ball pair
\[
(D^{n+1}, S^{n}, S^{n}\cap \R^{n+1}_+,S^{n}\cap \R^{n+1}_-).
\]
We can assume that the basepoint of $D^{n+1}$ is 
in $P^{n}\cap Q^{n} = S^{n-1}$.
For every object $W$ in $\V$  and $n\geq 0$ 
a ball pair in $\Vtop$ is defined as 
\[
(\square^{n+1}_W,S^n_W,P^{n}_W,Q^n_W) = 
W\times (\square^{n+1},S^n,P^{n},Q^n),
\]
and denoted simply by $(\square_W^{n+1},Q_W^n)$,
where $\square_W^{n+1} = D_{W}^{n+1}$.
If $Z$ is a finite $\V$-set, then 
we can define $\square^{n+1}_Z$ as the coproduct in $\Vtop$ of 
the $\square^{n+1}_W$ with $W\in Z$.
The pair
$(\square^{n+1}_Z,Q^{n}_Z)$ is also termed
\emph{ball pair} in $\Vtop$,
and $P^{n}_Z$ is termed the \emph{complement of $Q^{n}_Z$ 
in the boundary}.
It is easy to see that $\square^{n+1}_Z$ is a $\V$-complex
and $P^{n}_Z$ and $Q^{n}_Z$ are subcomplexes.
There are $\V$-cells only in dimensions $0$, $n-1$, $n$ and $n+1$:
the $0$-skeleton is $Z$, the $(n-1)$-skeleton
is $\partial Q^{n}_Z = \partial P^{n}_Z =  S^{n-1}_Z$,
the $n$-skeleton is $\partial \square^{n+1}_Z = S^n_Z  = 
P^{n}_Z\cup_{S^{n-1}_Z} Q^{n}_Z$, and the $(n+1)$-skeleton
is $\square^{n+1}_Z$ itself.

% later, in the proofs, use the fact that it is a cellular $I$-category.

Now consider a ball pair $(\square^{n+1}_Z,Q^{n}_Z)$
with the complement $P^{n}_Z$ in the boundary,
a $\V$-complex $L$ and a $\V$-map $f\from Q^{n}_Z \to L$
with the property that $f(Q^{n}_Z, S^{n-1}_Z) \subset 
(L_n,L_{n-1})$.
Let $K$ be the push-out in $\Vtop$
\ndiag{%
P^{n}_Z \ar[r]^f \ar@{ >->}[d] & L \ar@{ >->}[d]^{i} \\
\square^{n+1}_Z \ar[r] & K \\
}%
It is easy to check that $K$ is a $\V$-complex with subcomplex
$L$.
As in \cite{baues}, VIII.6, page 329, if $K'$ is a complex with $L$ as 
subcomplex
and $\V$-isomorphic to $K$ under $L$, then 
$K'$ is termed an \emph{elementary expansion} of $L$; moreover,
there is an associated  canonical $\V$-retraction $r\from K'\to L$ 
termed \emph{elementary collapse}.
If $S$ is a subcomplex of $L$ then we say that $i$ (resp. $r$)
is an elementary expansion (resp. collapse)
\emph{relative} $S$.

Now let $K$ and $L$ two finite $\V$-complexes
both containing a subcomplex $S$. 
If $j\from L \to K$ is an \emph{expansion relative $S$} 
% (see \cite{baues}, page 330)
we write $L\expansion K \rel S$;
is $r\from K \to L$ is a \emph{collapse relative $S$},
we write $K \collapse L \rel S$.
A \emph{formal deformation relative $S$}
(i.e. a finite composition of expansions and collapses
relative $S$)
is denoted by $L\deformation L$.
A $\V$-map $f\from L \to K$ under $S$ 
is a \emph{simple homotopy equivalence $\rel S$} 
if $f$ is $\V$-homotopic $\rel S$  to a formal deformation 
relative $S$.

Let $L$ be a finite $\V$-complex. 
Consider the set of all  pairs $(K,L)$ such that $L$ is a subcomplex
of $K$, $K$ is a finite $\V$-complex, and the inclusion $L\subset K$
is a $\V$-homotopy equivalence.
Two pairs $(K,L)$ and $(K',L)$ are equivalent if and only
if $K\deformation K' \rel L$. Let $[K,L]$ denote the 
class of the pair $(K,L)$ and 
let $\Wh(L)$ denote the set of equivalence classes.
It is possible to define an addition in $\Wh(L)$ by
\[
[K,L] + [K',L] = [K \cup_L K',L ]
\]
so that  $(\Wh(L),+)$ 
is an abelian group (lemma VIII.8.7, page 334 of \cite{baues}).
Given a cellular $\V$-map $f\from L \to  L'$ between 
finite $\V$-complexes, there is an induced map
\[
f_*\from \Wh(L) \to \Wh(L')
\]
defined by
\[
f_*[K,L] = [ K \cup_L M_f,L'],
\]
where $M_f$ is the mapping cylinder of $f$ in $\Vtop$ 
(see \ref{coro:factorization} below).
Actually $\Wh$ is a functor %$\Vcell_\equiv \to \Ab$ 
from the category of finite $\V$-complexes and $\V$-homotopy
classes of maps to the category of abelian groups.

Now consider a $\V$-homotopy equivalence $f\from X \to L$ 
of finite $\V$-complexes; the \emph{Whitehead torsion}
of $f$  is defined by 
\[
\torsion(f) = f_*[M_f,X] = [M_f\cup_X L,L]\in \Wh(L).
\]

\begin{theo}
\label{theo:torsion}
A $\V$-homotopy equivalence $f\from X \to L$ between finite 
$\V$-complexes is a simple homotopy equivalence if and only if 
the torsion vanishes $\torsion(f) = 0 \in \Wh(L)$.
Moreover:
\begin{enumerate}
\item For every pair of homotopy equivalences $f$ and $g$
the derivation property holds:
\[
\torsion(gf) = \torsion(g) + g_*\torsion(f).
\]
\item
Consider the following double push-out diagram
\ndiag{%
K_0' \ar[rrr] \ar[ddd] & & & K_2' \ar[ddd] \\
& K_0 \ar[ul]^{\equiv}_{f_0} \ar[r]^{i_2} \ar[d]^{i_1} \ar[rd]^{i_0} 
& K_2 \ar[ur]^{f_2}_{\equiv} \ar[d] \\
& K_1 \ar[dl]^{\equiv}_{f_1} \ar[r]  & K \ar@{.>}[dr]^{f}_{\equiv} \\
K_1' \ar[rrr] & & & K' \\ 
}%
where $f_0$, $f_1$ and $f_2$ are $\V$-homotopy equivalences
and all the non-diagonal maps are inclusions of subcomplexes.
Then the push-out map $f$ is a $\V$-homotopy equivalence %by \ref{lemma:ghost}
and the {addition formula} holds:
\[
\torsion(f) = \torsion(f_1\cup_{f_0} f_2) =
i_{1} \torsion(f_1) + i_2 \torsion(f_2) - i_0 \torsion(f_0).
\]
\end{enumerate}
\end{theo}
\begin{proof}
Let $\mathcal{D}$ denote the set of finite
$\V$-sets. 
Since $(\Vtop,\mathcal{D})$ is a cellular $I$-category (see \ref{rem:pro}),
theorem \ref{theo:torsion} is a consequence 
of VIII.8.3, VIII.8.4 and  VIII.8.5 (page 335) of \cite{baues}.
\end{proof}

Following Ranicki \cite{ranicki}, 
if $\AAA$ is a small additive category with sum denoted by $\oplus$,
then the \emph{isomorphism torsion
group} $\Kiso(\AAA)$
is the abelian group with one generator $\torsion(f)$ 
for each isomorphism $f\from M \to N$ in $\AAA$, and relations
\[
\torsion(gf)=\torsion(f) + \torsion(g)
\]
\[
\torsion(f\oplus f') = \torsion (f) + \torsion(f')
\]
for all isomorphisms 
$f\from M\to N$, $f'\from M' \to N'$ and $g\from N\to P$
in $\AAA$.
If $R$ is a ringoid, then it is possible to define 
$\Kiso(R)$ as the isomorphism torsion group of 
the additive category $\AAA(R)$ consisting 
of finitely generated free $R$-modules.
%Hence, by taking $R=\ze\Pi_\V(K)$, the group
%$\Kiso(R)$ is the isomorphism torsion group of the 
%category of the finitely generated free  $\ze\Pi_\V(K)$ modules.
%Consider in $\Kiso(R)$ the subgroup $T$ generated by 
%the invertible morphisms in $\Pi_\V(K)$,
%termed \emph{trivial units}.
A \emph{trivial} unit in $\Kiso(\ze\Pi_\V(K))$ is represented by 
an automorphism
of the free $\ze\Pi_\V(K)$-module $\ze\Pi_\V(K)(-,(V,x))$
that can be written as  $\pm f_*$ for some 
automorphism of $(V,x)$ in $\Pi_\V(K)$.

\begin{theo}
Let $K$ be a finite $\V$-complex. Then there is an isomorphism
\[
\tau\from \Wh(K) \cong 
K_1^{\mathbf{iso}}(\ze\Pi_\V(K))/{T}
%\Wh( \ze\Pi_\V(K)) 
\]
between the Whitehead group of $K$
and the quotient of $\Kiso(\ze\Pi_\V(K))$
by the subgroup $T$ generated by all trivial units.
It is defined by associating to $[X,K] \in 
\Wh(K)$ the torsion of the contractible 
$\Pi_\V(K)$-chain complex $C_*(X',K)$, 
where $X'$ is a normalization of $X$ rel $K$.
\end{theo}
\begin{proof}
It is a consequence of theorem VIII.12.7, page 348, of \cite{baues}.
\end{proof}

In case of transformation groups this is theorem 
14.16, page 286, of \cite{lueck} (in the same book the 
equivariant finiteness obstruction can be found).
Compare with the equivariant Whitehead torsion 
of Hauschild \cite{hauschild}, 
Dovermann and Rothenberg \cite{dovermannrothenberg},
Illman \cite{illman74} and Anderson \cite{anderson}.

%======================================================================
\section*{Appendix A. Basic homotopy theory in $\Vtop$}
\setcounter{equation}{0}
\renewcommand{\theequation}{(A.\arabic{equation}\/)}

We proved several results above by refering to 
\cite{baues}. 
We now describe some results that
are needed in order 
to apply the abstract theory of \cite{baues} to $\Vtop$.

\begin{defi}
\label{defi:cofibration}
A $\V$-map $i\from  Y \to X$ is a $\V$-cofibration if $i$ is a closed
inclusion and the following homotopy extension property holds:
for every commutative diagram in $\Vtop$
\ndiag{%
Y \ar@{ >->}[r]^i \ar@{ >->}[d]^{i_0} &
X \ar@{ >->}[d]^{i_0'} \ar@/^/[ddr]^{f} \\
Y\times I
\ar@{ >->}[r]^{i\times 1_I} \ar@/_/[rrd]_{H} &
X\times I
\ar@{.>}[dr]_G
\\
& &
Z
}%
there is a $\V$-map $G\from X \times I \to Z$ such that $G\circ i_0'=f$ and
$G\circ (i\times 1_I)=H$; in the diagram $i_0$ is the map
defined by $y \mapsto (y,0)$ for every $y\in Y$,
$i_0'$ the same for $X$.
\end{defi}

Equivalently a closed inclusion $i\from Y \to X$ in $\Vtop$
is a $\V$-cofibration if and only if the $\V$-map
$(Y\times I)\cup_Y X \to X\times I$  admits a retraction
(that is a left inverse). This implies that $\bar Y \to
\bar X$ and $Y \to X$ are cofibrations in $\top$.

\begin{lemma}\label{lemma:a2}
For a cofibration $i\from A\to M$ in $\Vtop$
and a $\V$-map $h\from A \to Y$,
in the push-out diagram
(which exists by \ref{lemma:pushout}) 
\ndiag{
A \ar[r]^h \ar[d]^i & Y \ar[d]^{\bar i} \\
M \ar[r] & X=M\cup_A Y \\
}%
the induced $\V$-map $\bar i$ is also a cofibration.
Moreover, $I$ carries this push-out diagram into a push-out diagram
(i.e. $(M\cup_A Y)\times I = (M\times I)\cup_{A\times I} (Y\times I)$).
\end{lemma}
\begin{proof}
The proof 
is formally identical to the proof for compactly generated spaces
with no structure group
(see (5.1) and (5.4) of \cite{whitehead}, pages 22--23);
or, directly applying the homotopy extension property \ref{defi:cofibration}
and the assumptions on the fibres, 
is the same as the proof in $\top$ (see e.g. lemma 2.3.6, page 56
of \cite{renzo}).
\end{proof}

The next result shows that basic homotopy theory is 
available in $\Vtop$.
For the axioms and properties of $I$-categories and cofibration
categories we refer the reader to \cite{baues}.
As in section 1 the cylinder in $\Vtop$ is given by
$IX = X\times I$.

\begin{theo}
The category $\Vtop$ with the cylinder $IX$
is an $I$-category.
Moreover, it is a cofibration category, with cofibrations
the $\V$-cofibrations and weak equivalences the $\V$-homotopy equivalences.
All objects are fibrant and cofibrant in $\Vtop$.
\end{theo}
\begin{proof}
To show that the axioms of
$I$-category are fulfilled, 
since lemma \ref{lemma:a2} holds, 
it suffices to use the argument of \cite{baues}, proposition (8.2),
with the same homeomorphism $\alpha\from I^2 \to I^2$, to prove the 
relative cylinder axion, see page 222 of \cite{baues}.
Then we can apply
Theorem (7.4) of \cite{baues}, page 223.
\end{proof}

The theorem implies that 
the following lemmata are true
in $\Vtop$. Corollaries \ref{coro:composition},
\ref{coro:apushout} and
\ref{coro:factorization} are immediate consequences of the first three
axioms of a cofibration category.

\begin{coro}[Composition]
\label{coro:composition}
Isomorphisms in $\Vtop$ are $\V$-homotopy equivalences and
also cofibrations.
For two maps
\ndiag{%
A \ar[r]^f & B \ar[r]^g & C \\
}%
if any two of $f$, $g$ and $gf$ are $\V$-homotopy equivalences
then so is the third. The composite of $\V$-cofibrations
is a $\V$-cofibration.
\end{coro}

\begin{coro}[Pushout]
\label{coro:apushout}
Let $i\from A\to M$ be a $\V$-cofibration and
$f\from A \to Y$ be a $\V$-map.
Consider the pushout in $\Vtop$,
\ndiag{%
A \ar@{ >->}[d]^i \ar[r]^f & Y \ar[d]^{\bar i} \\
M \ar[r]_{\bar f} & Y\cup_f M \\
}%
where 
$\bar i$ is a $\V$-cofibration by
\ref{lemma:a2}.
If $f$ is a $\V$-homotopy equivalence,
then so is $\bar f$;
if $i$ is a $\V$-homotopy equivalence, then so is $\bar i$.
\end{coro}

\begin{coro}[Factorization]
\label{coro:factorization}
Every $\V$-map $f\from X \to Y$
can be factorized as
\ndiag{%
X \ar[rr]^f \ar@{ >->}[dr]_i & & Y \\
&
Y' \ar[ur]^{\equiv}_g
}%
where $g$ is a $\V$-homotopy equivalence and $i$
a $\V$-cofibration.
The $\V$-family $Y'$ is 
the mapping cylinder (in $\Vtop$) of $f$, i.e. 
the push-out of $i_0\from X \to IX$ and $f$.%HERE
\end{coro}

\begin{coro}[Gluing lemma]
\label{coro:gluing}
Consider the following diagram,
\ndiag{%
A' \ar[rrr]^{h'} \ar@{ >->}[ddd] & & & Y' \ar[ddd] \\
& A \ar[ul]^{\equiv} \ar[r]^h \ar@{ >->}[d] & Y \ar[ur]^{\equiv} \ar[d] \\
& M \ar[dl]^{\equiv} \ar[r] & Y\cup_AM \ar@{.>}[dr]^{\equiv}  \\
M' \ar[rrr] &&& Y'\cup_{A'}M' \\
}%
where the two squares are push-outs in $\Vtop$, and
$A\to M$ and $A'\to M'$ are $\V$-cofibrations. The dotted
map exists as a consequence of \ref{lemma:a2}.
If the arrows $A\to A'$, $M\to M'$ and $Y\to Y'$ are $\V$-homotopy
equivalences, then so is the push-out map $Y\cup_AM \to Y'\cup_{A'}M'$.
\end{coro}
\begin{proof}
See lemma (1.2) of \cite{baues0}, page 84.
\end{proof}

\begin{coro}[Lifting lemma]
\label{coro:liftinglemma}
Consider a $\V$-cofibration $i\from A \to M$, a $\V$-homotopy
equivalence $p\from X\to Y$, and two maps $f\from A\to X$,
$g\from M\to Y$ such that $gi=pf$.
Then there is a $\V$-map $h\from M\to X$ such that
$hi=f$ and $ph\sim g \rel A$. The map $h$ is termed
the \emph{lifting} of the following diagram, and is unique up
to homotopy $\rel A$.
\ndiag{%
A\ar[r]^f \ar@{ >->}[d]^i & X \ar[d]^p_{\equiv}  \\
M\ar@{.>}[ur]^h \ar[r]^g & X \\
}%
\end{coro}
\begin{proof}
See lemma (1.11) of \cite{baues0}, page 90.
\end{proof}

The following is Dold's theorem in $\Vtop$.
\begin{coro}
Consider the commutative diagram.
\ndiag{%
X \ar[rr]^g&& Y \\
& A\ar@{ >->}[ul] \ar@{ >->}[ur] \\
}%
If $g$ is a $\V$-homotopy equivalence,
then $g$ is a $\V$-homotopy equivalence under $A$, that is,
there is a $\V$-map $f\from Y\to X$
such that $gf\sim 1_Y \rel A$
and $fg \sim 1_X \rel A$.
\end{coro}
\begin{proof}
See corollary (2.12) of \cite{baues0}, page 96.
\end{proof}

\begin{remark}
% \begin{lemma}\label{lemma:cofib}
Let $\V$ be a structure category with fibre functor $\ff$ and
$i\from A \to X$ a $(\V,\ff)$-cofibration.
Then $\bar i\from \bar A \to \bar X$ and $i\from A \to X$ are
(closed) cofibrations in $\top$. Moreover, since the fibre functor
$\ff$ is faithful, the induced map $i^V$ on function spaces
\[
i^V\from A^V \to X^V
\]
is a $(\V,\ff^V)$-cofibration for every $V\in \V$.
% \end{lemma}
\begin{proof}
Since
$i\from A \to X$
is a $(\V,\ff)$-cofibration, there is a
$(\V,\ff)$-retraction
\[
r\from X\times I \to A\times I \cup X\times\{0\}.
\]
By definition of $\V$-map, this implies that
$i$ and $\bar i$ are cofibrations in $\top$,
where clearly $A\subset X$ and $\bar A \subset \bar X$
are closed.
Moreover, assume that $\ff$ is faithful, so that the
function spaces $A^V$, $X^V$ are defined for every $V\in \V$.
Then there exists the  induced map
\[
r^V\from (X\times I)^V \to (A\times I \cup X\times\{0\})^V.
\]
Since $(X\times I)^V \cong X^V \times I$
and
\[
(A\times I \cup X\times\{0\})^V \cong
A^V\times I \cup X^V\times\{0\},
\]
and the homeomorphisms are actually
$(\V,\ff^V)$-isomorphisms,
there is
a $(\V,\ff^V)$-retraction
\[
X^V \times I \to
A^V\times I \cup X^V\times\{0\}.
\]
Since $i^V\from A^V \to X^V$ is a closed inclusion,
it  is a $(\V,\ff^V)$-cofibration.
\end{proof}
\end{remark}

% definition of strong deformation retract
We recall now the definition of $\V$-deformation retract,
following Whitehead \cite{whitehead}, page 24.
A $\V$-family $A$ is a $\V$-deformation retract of $X\supset A$
if there is a $\V$-homotopy $F\from X\times I \to X$
such that
$F(x,t)=x$ for each $(x,t)\in A\times I \cup X\times \{0\}$
and
$F(X\times\{1\}) \subset A$.
%This is what some authors called \emph{strong}
%deformation retract, if $\V$ is trivial.

\begin{lemma}\label{lemma:ndr}
Let $X$ be a $(\V,\ff)$-family and
$A\subset X$ a closed inclusion.
Then for every $V\in \V$
the induced map $A^V \to X^V$ is a closed inclusion
of $(\V,\ff^V)$-families.
Moreover, if $A$ is a $(\V,\ff)$-deformation retract of $X$,
then $A^V$ is a $(\V,\ff^V)$-deformation retract of $X^V$.
$A$ is a deformation retract of $X$ in $\top$
and $\bar A$ is a deformation retract of $\bar X$ in $\top$.
\end{lemma}
\begin{proof}
The inclusion $A^V \to X^V$ is a closed inclusion
since $A^V\subset X^V$ is equal to the pre-image
of $\bar A$
under the map $X^V \to \bar X$.
Furthermore, any
$(\V,\ff)$-deformation retraction $F\from X\times I \to X$
yields a $(\V,\ff^V)$-deformation retraction
$F^V\from X^V\times I \to X^V$.
\end{proof}

\setcounter{equation}{0}
\renewcommand{\theequation}{(B.\arabic{equation}\/)}
\section*{Appendix B. The category of coefficients}
We here describe the coefficient functor for $\Vtop$ which 
is a basic ingredient of the theory in \cite{baues}.

\begin{defi}
\label{defi:T}
Let $A$ be a $\V$-space. Consider pairs of $\V$-maps 
$x,y\from Z \to A$
where $Z$ is a $\V$-set. We associate with $x, y$
the $\V$-space $X(x,y)$ 
given by the push-out in $\Vtop$
\ndiag{%
Z\times \{0,1\} \ar[r]^{\hspace{24pt}(x,y)} \ar[d] & A \ar[d] \\
Z\times I \ar[r] & X(x,y) \\
}%
\end{defi}
Hence $(X(x,y),A)$ is a reduced $1$-dimensional relative
$\V$-complex.
We call it a \emph{$\V$-graph under $A$}.
Let $\TT$ denote the category \emph{of $\V$-graphs}:
Objects are $\V$-graphs $X=(X,A) = (X(x,y),A)$ 
and morphisms are $\V$-homotopy classes rel $A$ of $\V$-maps,
so that 
\[
\TT(X,Y) = [X,Y]^A_\V
\]
for $\V$-graphs $X$, $Y$ under $A$. The coproduct
$X\vee_A Y = X\cup_A Y$ in $\TT$ is defined by the push-out
of $X\leftarrow A \to Y$ in $\Vtop$.
Since finite coproducts exists in $\TT$, $\TT$ is a theory.
Let $\TT^\sharp$ denote the full subcategory of $\TT$ 
consisting of $\V$-graphs $X(x,y)$ for which the 
$\V$-set $Z$ above is finite. Since $\V$ is a small category,
we see that $\TT^\sharp$ is a small category.
Let $\model(\TT^\sharp)$ denote the \emph{category 
of models of $\V$-graphs}, that is the category of functors
$\left(\TT^{\sharp}\right)^\op \to \Set$
which carry coproducts in $\TT^\sharp$
to products in $\Set$.

\begin{lemma}
\label{lemma:Tcomplex}
A $\TT$-complex in $\Vtop^A$ 
as defined in \cite{baues}
is the same as a relative $\V$-complex
$(X,A)$ which is reduced and normalized.
\end{lemma}
\begin{proof}
According to IV.2.2 of \cite{baues} a $\TT$-complex
$X_{\geq 1}$ in $\Vtop$ is given by a sequence of $\V$-cofibrations
\[
X_{\geq 1} = (X_1 \subset X_2 \subset \dots)
\]
where $X_1$ is an object in $\TT$ and $(X_{n+1},X_n)$ 
is a principal cofibration with attaching map
$\partial_{n+1}\from \Sigma^{n-1} C_{n+1} \to X_n$, $n\geq 1$,
where $C_{n+1}$ is a cogroup in $\TT$
(that is, a $1$-dimensional spherical object $S^1_\alpha$
under $A$). Since $\Sigma^{n-1} S^1_\alpha = S^n_\alpha$,
we have the result.
\end{proof}

%\begin{lemma}
\begin{remark}
\label{rem:pro}
The category $\Vtop$ is a cofibration category under $\mathbf{T}$
that satisfies the \emph{Blakers--Massey property}
%\end{lemma}
%\begin{proof}
(see \cite{baues}, page 245).
%\end{proof}
%\begin{lemma}
%\label{lemma:homological}
Hence
it follows from proposition 1.2, page 250 of \cite{baues}
that
$\Vtop$ is a homological cofibration category under $\mathbf{T}$.
Moreover 
%\label{lemma:cellularcategory}
$\Vtop$ is a cellular $I$-category since 
the ball pair axiom holds, as a consequence of the ball pair axiom in
$\top$  and 
the cellular approximation theorem 
\ref{theo:cat}
(see page 327 of \cite{baues}).
\end{remark}

% track category
We recall that a \emph{track category},
i.e. a \emph{category enriched in groupoids},
is a $2$-category all of whose $2$-cells  are invertible
(see \cite{kelly} for details).
In particular the category of groupoids $\grd$  is a track category.
Objects are groupoids, morphisms are functors between groupoids
and tracks are natural transformations.
A \emph{track functor}, or else a $2$-functor, 
between track categories is a groupoid enriched functor.

Let $\V$ be a topological enriched category. Then there is 
a canonical track category $\llbracket \V \rrbracket$ associated to $\V$
defined as follows. Objects and morphisms are the same as $\V$
and $2$-cells are tracks 
$b\from \alpha \tto \alpha'$ in $W^V = \hom_\V(V,W)$
with 
$\alpha,\alpha'\from V \to W$ and $V$, $W$ objects in $\V$

A $\llbracket \V^\op \rrbracket$-diagram 
is a track functor $\llbracket \V^\op \rrbracket \to \grd$.
A morphism between $\llbracket \V^\op \rrbracket$-diagrams is a
natural transformation of track functors.
An example of a track functor is given by
\[
\llbracket \Pi_\V(X^\circ,A^\circ) \rrbracket 
\from \PiV \to \grd,
\]
with $(X,A)$ a pair of $\V$-families.
The image of an object $V\in \PiV$ is the groupoid 
$\Pi(X^V,A^V)$. The image of a morphism $\lambda\from V \to W$ in 
$\PiV$ is the morphism of groupoids
$\lambda^*\from \Pi(X^W,A^W) \to \Pi(X^V,A^V)$
induced by the map $X^\lambda\from X^W \to X^V$.
The image of a track $b\from \lambda \tto \lambda'$ in 
$\PiV$ is given by the natural transformation
\ndiag{%
\Pi(X^W,A^W) \ar@/^2pc/[rr]^{\Pi(\lambda)} \ar@/_2pc/[rr]_{\Pi(\lambda')}  &
\ar@{} |{\displaystyle\Downarrow b^*} & 
\Pi(X^V,A^V)\\
}%
in $\gpd$.
%\begin{comment}
%More explicitely, if $x_W,y_W\from W\to A$ are $\V$-points in $A$,
%$V$ and $W$ are objects in $\V$, $\lambda$ and $\lambda'\from V \to W$
%are morphisms in $\V$ and $b\from \lambda \tto \lambda'$ is a track,
%then the induced track 
%associate to $x_W$ the track
%$x_Wb \in \Pi(X^V,A^V)(x_W\lambda,x_W\lambda')$
%and to $y_W$ the track
%$y_Wb\in \Pi(X^V,A^V)(y_W\lambda,y_W\lambda')$.
%There is also a morphism
%\[
%b^*\from \Pi(X^V,A^V)(x_W,y_W) \to \Pi(X^W,A^W)(x\lambda,y\lambda')
%\]
%given by the composition
%\[
%(\xi\from x \tto y)
%\mapsto
%(\xi b \from x\lambda\tto y\lambda').
%\]
%% todo: it is the composition of the squares in the natural transformation!
%For every track $a\from x_W \tto y_W$ the image  $b^* a$ 
%is equal to the diagonal of the following commutative square.
%\ndiag{%
%x_W\lambda \ar@{=>}[r]^{x_Wb} \ar@{=>}[d]_{a\lambda}  
%\ar@{=>}[rd]^{b^*a}
%& x_W\lambda' \ar@{=>}[d]^{a\lambda'} \\
%y_W\lambda \ar@{=>}[r]_{y_Wb} & 
%y_W\lambda' \\
%}%
%\end{comment}

With the aid of the track category $\PiV$ we define the category 
$\FAgrd$ as follows.
Objects are track functors 
$F\from \llbracket\V^\op\rrbracket \to \grd$ 
under $\llbracket \Pi(A^\circ) \rrbracket $
such that the natural transormation
\[
\iota\from \llbracket \Pi(A^\circ) \rrbracket \to F
\]
for each $V\in \V$ yields a functor
\[
\iota \from \Pi(A^V) \to F^V
\]
which is the identity on objects.
A morphism in $\FAgrd$ is a natural transformation 
of track functors under $\llbracket \Pi(A^\circ) \rrbracket $.

\begin{theo}
\label{theo:ghost}
Let $\TT^\sharp$ denote the  theory of finite $\V$-graphs under $A$. Then
there is an equivalence of categories
\[
\model(\TT^{\sharp}) \equiv \FAgrd.
\]
\end{theo}
\begin{proof}
Consider the following functor
\[
\xi\from \model(\TT^{\sharp}) \to \FAgrd.
\]
Let $M\from (\TT^\sharp)^\op \to \Set$ be a model
of $\TT^\sharp$.
Then $\xi(M)$ is the $\llbracket \V^\op \rrbracket$-diagram
of groupoids under $\llbracket \Pi(A^\circ) \rrbracket$ 
defined as follows. Let $V$ be an object of $\V$.
Objects 
of the groupoid $\xi(M)(V)$ 
are the same objects as $\Pi_\V^V(A) = \Pi(A^V)$, i.e. 
$\V$-points of type $x_V\from V\to A$.
For every pair of objects $x_V,y_V$  in $\xi(M)(V)$ 
the homset is given by the set $M(X(x_V,y_V))$.
Let $x_V,y_V$ and $z_V$ be points in $A^V$. 
The natural $\V$-map
\[
X(x_V,z_V) \to X(x_V\cup y_V, y_V \cup z_V) = X(x_V,y_V) \vee X(y_V,z_V)
\]
induces under $M$ a function
\ndiag{%
M(X(x_V\cup y_V, y_V \cup z_V)) 
\ar@{=}[d] \\
 M(X(x_V,y_V)) \times M(X(y_V,z_V)) \ar[r] & 
 M(X(x_V,z_V))\\
}%
which yields the 
composition law 
for morphisms in $\xi(M)$. 
Since there is a natural $\V$-isomorphism $X(x_V,y_V) \equiv X(y_V,x_V)$,
we obtain that $\xi(M)(V)$ is a groupoid.

Now let $\lambda\from V\to W$ be a morphism in $\V$.
The induced $\V$-map $A^W \to A^V$ is a function on objects
of $\xi(M)$. Given $x_W$ and $y_W$ in $A^W$ 
the induced $\V$-map 
\[
X(x_W\lambda,y_W\lambda) \to X(x_W,y_W)
\]
has an image under $M$ 
\ndiag{%
\xi(M)(W)(x_W,y_W) \ar@{=}[d] &  
 \xi(M)(V)(x_W\lambda,y_W\lambda) \ar@{=}[d] \\ 
M(X(x_W,y_W)) \ar[r]   & 
M(X(x_W\lambda,y_W\lambda))
\\
}%
Hence $\xi(M)(\lambda)$ is a morphism of groupoids,
and $\xi(M)$ is a functor $\V^\op \to \gpd$. It is not difficult 
to show that it is a functor under $\Pi(A^\circ)$.

It is left to define the values of $\xi(M)$  on the tracks 
of $\llbracket \V^\op \rrbracket$.
A track $b\from \lambda \tto \lambda'$ 
in $\llbracket \V^\op \rrbracket$ is
represented by a $\V$-map
\[
b\from V\times I \to W.
\]
Hence for every $x_W \from W \to  A$ the track $b$ 
induces a $\V$-map $x_Wb\from V\times I \to A$, and hence 
a $\V$-map
\[
X(x_W\lambda,x_W\lambda') \to A.
\]
By taking the values of the model $M$ we obtain 
a function
\[
* = M(A)  \to 
M(X(x_W\lambda,x_W\lambda')),
\]
hence we obtain an element in $M(X(x_W\lambda,x_W\lambda'))$ which 
is therefore a morphism in 
$\xi(M)(V)(x_W\lambda,x_W\lambda')$.
This defines the  natural transormation induced by $b$. 
By the elementary properties of $M$ it is easy to prove that 
$\xi$ is a track functor.

Conversely, let $\CC$ be an object in $\FAgrd$. 
Then to each $\V$-graph $X(x_Z,y_Z)$
with $x_Z,y_Z\from Z\to A$ we can associate the set
\[
\xi'(\CC)(X(x_Z,y_Z) = 
\Mor\left( 
\llbracket \Pi(X(x_Z,y_Z)^\circ,A^\circ) \rrbracket , \CC
\right)
\]
of all morphisms in $\FAgrd$ from 
$\llbracket \Pi(X(x_Z,y_Z)^\circ,A^\circ) \rrbracket $ to $\CC$.
It is clear that $\xi'(\CC)$ is a functor
$\left( \TT^\sharp \right)^\op \to \Set$. 
It is also a model by \ref{propo:vankampen} below.
% being a push-out it is possible to ...
It is possible to prove that this construction
yields a functor
\[
\xi'\from 
\FAgrd
\to 
\model(\TT^{\sharp}).
\]
The functors $\xi$ and $\xi'$ are equivalences
of categories (the details are left to the reader).
\end{proof}

We define the 
\emph{coefficient functors} $c$ and $\pi_\V$  
\ldiag{diag:coeff}{%
%\pi_\V \from 
& \model(\TT^\sharp) \ar[d]^{\equiv}_{\xi}\\
(\Vtop)^A/_{\equiv \rel A} \ar[r]^{\hspace{12pt}\pi_\V} \ar[ur]^{c} &  
\FAgrd \\
}%
as follows.
For every $\V$-family $(X,A)$ under $A$ (where 
$A\to X$ is a $\V$-cofibration)
let $\pi_\V(X,A)  = \llbracket \Pi(X^\circ,A^\circ) \rrbracket $.
Furthermore, let $c(X,A)$ be the model in $\model(\T^\sharp)$ 
that associates to  the $\V$-graph
$X(x,y)$ with $x,y\from Z \to A$ the set of $\V$-homotopy classes $\rel A$ 
of $\V$-maps extending $1_A$ 
\[
X(x,y) \mapsto [X(x,y),X]^A_\V.
\]

% define c

\begin{propo} The diagram \ref{diag:coeff} is commutative:
\[
\xi c = \pi_\V.
\]
\end{propo}
\begin{proof}
The object
$\xi c(X,A)$ in $\FAgrd$ associates to $V \in \llbracket \V^\op \rrbracket$
the groupoid $\xi c(X,A)^V$ with objects the points in $A^V$; the morphisms
are the sets 
\[
c(X,A)(X(x,y)) = [X(x,y),X]^A_\V
\]
where $x,y\from V\to A$ are  $\V$-points.
The morphism set $[X(x,y),X]^A_\V$
is naturally isomorphic to $\Pi(X^V,A^V)(x,y)$. 
The same argument can be applied to tracks, so that 
we get $\xi c  = \pi_\V$ as claimed. 
\end{proof}

Now we define the category $\FAset$.
Objects are pairs of $\V$-maps $x,y\from Z \to A$
where $Z$ is a $\V$-set.
A morphism 
$(x,y) \to (x',y')$ 
is a $\V$-map $u\from Z \to Z'$ 
such that $x'u=x$ and $y'u=y$.

Consider now an object  $F\from \llbracket \V^\op \rrbracket \to 
\grd$ in $\FAgrd$. For every object $V\in \V$ we denote by $F^V$ the 
groupoid associated to $V$.
For every $V\in \V$ consider the two maps
\ndiag{
F^V_1 
\ar@<2pt>[r]^{\partial_0}  
\ar@<-2pt>[r]_{\partial_1}  
& 
A^V\\
}%
from the set of arrows of $F^V$ to $A^V$ given by the 
source $\partial_0$ and the target $\partial_1$ 
in the groupoid $F^V$.
This yields an object $x_V,y_V\from Z_V = V\times F^V_1 \to  A$ in $\FAset$.
The coproduct of all such objects is 
an object $\Phi(F)= x,y\from Z \to A$. 
This defines the forgetful functor
\[
\Phi\from \FAgrd \to \FAset.
\]
The left-adjoint of $\Phi$ 
denoted by $<,>$ sends an object $x,y\from Z\to A$
in $\FAset$ to the free object $<x,y>$ in $\FAgrd$.

\begin{propo}
For every $x,y\from Z \to A$, we have a natural isomorphism
\[
<x,y> = \pi_\V(X(x,y),A).
\]
\end{propo}

Moreover a result  corresponding to the Seifert--van Kampen theorem 
holds. The proof is 
a modification of the Brown's proof \cite{brown,brown2}
of the Seifert--van Kampen theorem for groupoids.

\begin{propo}
\label{propo:vankampen}
If $X_1$ and  $X_2$ are 
$\V$-complexes relative $A$ both containing 
a copy of a relative $\V$-complex $(X_0,A)$ as a subcomplex
then
\[
\pi_\V(X_1\cup_{X_0} X_2) = \pi_\V(X_1) \cup_{\pi_\V(X_0)} \pi_\V(X_2).
\]
Here the right hand side  denotes the push-out in $\FAgrd$.
\end{propo}

\begin{remark}
In \cite{baues} the coefficient functor $c$ maps to the category
$\coef$ which can be identified with
\[
\FAgrd \equiv \model(\TT^\sharp) \equiv \coef
\]
such that the coefficient functor $c$ in \cite{baues} 
mapping to $\coef$ coincides with $\pi_\V$ and $c$ in \ref{diag:coeff}.
Moreover 
the \emph{enveloping functor} $U$ of \cite{baues}
(pages $\geq 150$) %in case of $\Vtop$ 
can be identified with the composite
\[
U \from 
\xymatrix@1{\coef \ar[r]^{\equiv} & 
\FAgrd \ar[r]^{\bar U} & 
\ringoids
}
\]
where $\bar U$ carries $F$ to the category
$\ze\left[ (\int_\V F)/_\sim\right]$.
Here the equivalence relation $\sim$ for the integration category
$\int_\V F$ is defined as in \ref{eq:fundcat} above. Hence 
for a reduced relative $\V$-complex $(X,A)$ we get
\[
Uc(X,A) = \bar U (\pi_\V(X,A) ) = \ze\Pi_\V(X,A)
\]
where we use the discrete fundamental category in \ref{eq:fundcat}.
\end{remark}

%======================================================================

\def\cfudot#1{\ifmmode\setbox7\hbox{$\accent"5E#1$}\else
  \setbox7\hbox{\accent"5E#1}\penalty 10000\relax\fi\raise 1\ht7
  \hbox{\raise.1ex\hbox to 1\wd7{\hss.\hss}}\penalty 10000 \hskip-1\wd7\penalty
  10000\box7} \def\cprime{$'$}

%======================================================================

\end{document}